\documentclass[hidelinks,onefignum,onetabnum]{siamart220329}
\usepackage{lipsum}
\usepackage{amsfonts}
\usepackage{graphicx}
\usepackage{epstopdf}
\usepackage{algpseudocode}
\Crefname{ALC@unique}{Line}{Lines}
\ifpdf
  \DeclareGraphicsExtensions{.eps,.pdf,.png,.jpg}
\else
  \DeclareGraphicsExtensions{.eps}
\fi

\usepackage{amsopn}
\usepackage{amsmath}
\usepackage{amssymb}
\usepackage{mathrsfs}
\usepackage{hyperref}
\usepackage{bbm}
\usepackage{graphicx}
\usepackage{float}
\usepackage{booktabs}
\usepackage{color}
\usepackage{subfigure}
\usepackage{mathtools,tikz-cd} 
\usepackage{pst-plot}
\usepackage{algorithm}
\usepackage{parskip}
\newcommand{\ud}{\mathrm{d}}
\newcommand{\bd}{\boldsymbol}

\newcommand{\ee}{\ensuremath{\text{e}}}
\newcommand{\grad}{\mathrm{grad}\,}

\newcommand{\Tr}{\mathrm{Tr}\,}
\newcommand{\st}{\mathrm{s.t.}\,}

\newcommand{\bmat}[1]{\begin{bmatrix} #1 \end{bmatrix}}

\newcommand{\norm}[1]{\left\lVert#1\right\rVert}
\newcommand{\abs}[1]{\left \lvert #1 \right \rvert}

\newcommand{\E}[1]{\mathbbm{E} \left( #1 \right)}

\newcommand{\energy}{\mathcal{E}}
\newcommand{\Snp}{\mathcal{S}^{n,p}_+}
\newcommand{\R}{\mathbb{R}}
\newcommand{\rhoref}{\rho_{\mathrm{ref}}}
\newcommand{\Ito}{It\^{o}}

\newtheorem{thm}{Theorem}

\newtheorem{lem}[thm]{Lemma}
\newtheorem{rem}[thm]{Remark}
\newtheorem{prop}[thm]{Proposition}

\newcommand{\StatexIndent}[1][3]{%
  \setlength\@tempdima{\algorithmicindent}%
  \Statex\hskip\dimexpr#1\@tempdima\relax}
\newcommand{\appropto}{\mathrel{\vcenter{
  \offinterlineskip\halign{\hfil$##$\cr
    \propto\cr\noalign{\kern2pt}\sim\cr\noalign{\kern-2pt}}}}}
\graphicspath{{pics/}}

\newcommand*\samethanks[1][\value{footnote}]{\footnotemark[#1]}

\headers{Riemannian Langevin Monte Carlo for PSD matrices of fixed rank}{T. Yu,  S. Zheng, J. Lu, G. Menon and X. Zhang}

\title{Riemannian Langevin Monte Carlo schemes for sampling PSD matrices with fixed rank \thanks{\funding{S.Z. and X.Z. are supported by NSF DMS-2208518. J.L. is supported in part by NSF DMS-2012286 and DMS-2309378.  G.M. is supported in part by NSF DMS-2107205}}}

\author{Tianmin Yu\thanks{Division of Applied Mathematics, Brown University,    Providence, RI (\email{tianmin\_yu@brown.edu,  govind\_menon@brown.edu}).}
\and
Shixin Zheng\thanks{Department of Mathematics,
Purdue University,
West Lafayette, IN
(\email{zheng513@purdue.edu, zhan1966@purdue.edu}).}
\and Jianfeng Lu\thanks{Departments of Mathematics, Physics, and Chemistry, Duke University, Durham,
NC (\email{jianfeng@math.duke.edu}).}
\and Govind Menon\samethanks[2]
\and Xiangxiong Zhang\samethanks[3]
}

\begin{document}
\maketitle

\begin{abstract}
This paper introduces two explicit schemes to sample matrices from Gibbs distributions on $\mathcal S^{n,p}_+$, the manifold of real positive semi-definite (PSD) matrices of size $n\times n$ and rank $p$. Given an energy function $\mathcal E:\mathcal S^{n,p}_+\to \mathbb{R}$ and certain Riemannian metrics $g$ on $\mathcal S^{n,p}_+$, 
these schemes rely on an Euler-Maruyama discretization of the Riemannian Langevin equation (RLE)
with  Brownian motion on the manifold. 
We present numerical schemes for RLE under two fundamental metrics on $\mathcal S^{n,p}_+$: (a) the metric obtained from the embedding of $\mathcal S^{n,p}_+ \subset \mathbb{R}^{n\times n} $; and (b) the Bures-Wasserstein metric corresponding to quotient geometry. We also provide examples of energy functions with explicit Gibbs distributions that allow numerical validation of these schemes.
\end{abstract}


% REQUIRED

\section{Introduction}
\label{sec:intro}
\subsection{Problem statement}
Consider the space of real, symmetric positive semi-definite matrices with size $n\times n$ and rank $p$, denoted by
\begin{equation}
\label{eq:defSnp}
\mathcal{S}^{n,p}_+=\{X\in \mathbb R^{n\times n}|X=X^T,X\succeq 0,\mathrm{rank}(X)=p\}.    
\end{equation} 
Given an energy $\energy: \Snp \to \R$ and a parameter $\beta>0$ referred to as the inverse temperature, our goal is to sample efficiently from the Gibbs distribution
\begin{equation}
\label{eq:def-Gibbs}
    \rho_\beta(X) = \frac{1}{Z_\beta} e^{-\beta \energy(X)} \rhoref(X), \quad Z_\beta = \int_{\mathcal{S}^{n,p}_+} e^{-\beta \energy(X')} \rhoref(X') \, dX'.
\end{equation} 
Gibbs measures must be defined with respect to a base measure. In this work, we equip the space $\Snp$ with a Riemannian metric $g$ and choose $\rhoref(X) dX= \sqrt{\det g(X)} dX$ to be the canonical volume form associated  to the metric $g$. This volume form is expressed in coordinates for the metrics studied in this paper in Section~\ref{sec:example}.

This sampling problem is related to the optimization problem $\min_{X \in \Snp}\energy(X)$ since in the limit 
$\beta\to \infty$ the Gibbs distribution concentrates at the global minima of $\energy (X)$. Minimization problems over the space $\mathcal{S}^{n,p}_+$ arise in many areas, especially semidefinite programming and machine learning, and have been studied extensively. Gibbs distributions originate in statistical physics, while the sampling problem may also be seen as a stochastic variant of the optimization problem.
%$\min_{X \in \Snp}\energy(X)$. 
For these reasons, the sampling problem has a broad range of applications; see Section~\ref{subsec:applications} below.

The main contribution of this paper are efficient sampling schemes  
for $\rho_\beta$ based on Langevin dynamics. Our approach builds on the geometric theory of optimization; in particular, we extend Riemannian optimization on $\mathcal S_+^{n,p}$ ~\cite{vandereycken2009embedded,zheng2022riemannian} to Gibbs sampling as follows. In~\cite{vandereycken2009embedded} it was recognized that two commonly used gradient descent schemes over $\Snp$ are time discretizations of {\em Riemannian\/} gradient flows, where $\Snp$ is equipped with the two natural Riemannian metrics listed below. We combine this observation with the theory of Brownian motion on Riemannian manifolds to obtain Riemannian Langevin equations and explicit sampling schemes. 

The reader unfamiliar with these concepts should note that while the abstract theory serves to guide our work, the schemes presented in this paper may be implemented without requiring a complete understanding of the underlying theory. Further, while this paper is focused on the two numerical schemes below, the underlying framework can be used to extend other Riemannian gradient descent schemes to sampling schemes for the Gibbs measure. The new phenomenon that arises is the interplay between Brownian motion and curvature in the Riemannian Langevin equation. This interplay has been studied in depth by two of the authors (TY and GM) and their co-workers in recent papers for geometries used in optimization and physics~\cite{IM,MY1,MY2}.

%We resolve these issues for $(\Snp,g)$ for the two fundamental metrics $g$ arising in optimization over $\Snp$~\cite{massart2020quotient,massart2019curvature,vandereycken2009embedded}. 

\subsection{Two Riemannian metrics on $\mathcal{S}^{n,p}_+$}
Given $X\in \mathcal{S}^{n,p}_+$, let  $X=YY^T$ be a low-rank decomposition where $Y\in \mathbbm R^{n\times p}$. We use two fundamental metrics on $\Snp$ obtained from this parametrization, from the Euclidean metric for either the variable $X$ or the variable $Y$ through the use of Riemannian embedding and Riemannian submersion respectively. These are the two most natural ways of defining metrics on $\Snp$.

The flat metric for $X$ corresponds to the embedded geometry of  $\mathcal{S}^{n,p}_+$ in the Euclidean space $\mathbbm R^{n\times n}$~\cite{vandereycken2009embedded}. Precisely, we consider the natural Riemannian embedding $\mathcal{S}^{n,p}_+\xhookrightarrow{}  \mathbbm R^{n\times n}$ and use the Frobenius norm on $\mathbbm R^{n\times n}$ to define a metric on $\Snp$.  
Denote it by $g_E$, then $g_E(A,B)=\Tr(A^TB)$ for any two square matrices $A, B$ in the tangent space of $\mathcal{S}^{n,p}_+$, where $\mbox{Tr}$ denotes the trace of a matrix.

On the other hand, we may also use the flat geometry on $Y$ to define a metric on $\Snp$. We observe that if $YY^T=X$, then it is also true that $\tilde{Y}\tilde{Y}^T=X$ where $\tilde{Y}=YO$ and $O \in \mathcal{O}_p$, the orthogonal group of dimension $p$. Thus, we may identify $\mathcal{S}^{n,p}_+\eqsim \mathbb{R}^{n\times p}_*/\mathcal{O}_p,$ as a quotient space, with a quotient map \[\begin{aligned}
    \pi:  \mathbb{R}^{n\times p}_* &\rightarrow \mathbb{R}^{n\times p}_*/\mathcal{O}_p\\
     Y &\mapsto [Y]=\left\{ YO \mid O\in \mathcal{O}_p \right\}. 
\end{aligned}\]  
Here $\mathbb{R}^{n\times p}_*$ denotes full rank matrices.

The quotient space structure can be enhanced with a Riemannian metric through the use of Riemannian submersion. Roughly, the metric for $X$   corresponds to the metric for $Y$ in a manner that respects the splitting of the tangent space at $Y$ into the space of the group action and its complement. If $\mathbb{R}^{n\times p}_*$ is equipped with Euclidean metric, then the metric induced by the submersion is often called the  Bures-Wasserstein metric on $\mathcal{S}^{n,p}_+\eqsim \mathbb{R}^{n\times p}_*/\mathcal{O}_p$, denoted by $g_{BW}$~(see~\cite{bhatia2017bureswasserstein,massart2020quotient,massart2019curvature}).

\subsection{Langevin dynamics and the Riemannian Langevin equation}
We now explain how Langevin equations may be defined intrinsically on $(\Snp,g)$. 

Let us first recall the  Langevin equation on $\mathbb{R}^n$. Assume given a potential or energy function $\energy: \mathbb{R}^n\to\mathbb R$ and let $W_t$ denote the standard Wiener process on $\mathbb{R}^n$. The Langevin equation for the potential $\energy$ is the \Ito\/ differential equation
\begin{equation}
    \label{eq:langevin}
    dx_t = -\nabla \energy(x_t) \, dt + \sqrt{\frac{2}{\beta}} \, dW_t.
\end{equation}
The Fokker-Planck equation describes the evolution of the probability density of $x_t$. With $\rho(x,t)\, dx = \mathbb{P}(x_t \in (x,x+dx))$, we have
\begin{equation}
    \label{eq:fp}
    \partial_t \rho = \frac{1}{\beta} \triangle \rho + \nabla \cdot \left(\rho \nabla \energy\right).
\end{equation}
The Gibbs density (with reference density being uniform with respect to Lebesgue measure) is the unique equilibrium of equation~\eqref{eq:fp} under natural growth assumptions on the energy $\energy$ as $|x|\to \infty$.
% \jl{strictly speaking, this requires some growth assumption of $\energy$}

The Langevin equation immediately yields a numerical scheme for (approximate) sampling from the Gibbs distribution. Fix a step size $\Delta t >0$, let $t_k = k \Delta t$, $k=0,1,\ldots$, and let $x_k$ denote the numerical approximation to~\eqref{eq:langevin} at time $t_k$. The Euler-Maruyama scheme to approximate equation~\eqref{eq:langevin}, also known as Langevin Monte Carlo in the statistics literature, is
\begin{equation}
    x_{k+1}=x_k-\Delta t \nabla \mathcal E(x_k)+\sqrt{\frac{2\Delta t}{\beta}}\xi_k,
    \label{scheme-LD}
\end{equation}
where $\xi_k =(\xi_k^1, \ldots, \xi_k^n)$ is an i.i.d. sequence of standard Gaussian vectors in $\mathbb{R}^n$. This scheme is explicit. In order to extend it to sampling from~\eqref{eq:def-Gibbs} we must understand how to modify the Langevin equation on the Riemannian manifold $(\Snp,g)$. 

First, the term $\nabla \energy$ must be replaced by the Riemannian gradient, written as $\grad \energy$.
The more subtle modification of equation~\eqref{eq:langevin} concerns the noise. The natural analogy is to replace the Wiener process $W_t$ on $\mathbb{R}^n$ with Brownian motion on the Riemannian manifold $(\Snp,g)$ at inverse temperature $\beta$, denoted $\bd B_t^{g,\beta}$. This yields the (formal) Riemannian Langevin equation on $(\Snp,g)$ 
\begin{align}\label{sde1}
    \ud \bd X_t=-\grad \mathcal E(\bd X_t)\ud t+\ud \bd B_t^{g,\beta}.
\end{align}

%defined by
%\begin{equation}
%    g\left(\grad \energy, v\right) = \left. \frac{d}{d\tau} \energy \left(x(\tau)%\right)\right|_{\tau =0},
%\end{equation}
%where $x(\tau)$ is a smooth curve on $\Snp$ such that $\dot{x}(0)=v \in T_x \Snp$. %The metric `converts' the differential of $\energy$ into the gradient vector. 

%As noted previously, the role of the Riemannian gradient for optimization over $\Snp$ has been studied in~\cite{massart2019curvature,massart2020quotient,vandereycken2009embedded}. 

%Here $\bd X_t$ is a random process taking values in $\mathcal{S}^{n,p}_+$ and $\bd B_t^{g,\beta}$ stands for the Brownian motion on $(\mathcal{S}^{n,p}_+,g)$ at inverse temperature $\beta$. 

%The precise meaning of such a Brownian motion term shall be made clear in Section \ref{sec:SDE}
%\begin{equation}
%    \label{eq:RLE-abstract2} dX_t = -\grad_g E(X_t) + dB_t^{g,\beta},
%\end{equation}
%where $B_t^{g,\beta}$ denotes  Brownian motion on $(\mathcal S^{n,p}_+,g)$ at inverse %temperature $\beta$. 

This equation is only formal because stochastic differential equations on manifolds must be defined using the Stratonovich formulation in order to ensure coordinate independence (\Ito\/ differentials do not satisfy the chain rule, while Stratonovich differentials do)~\cite{Hsu,Ikeda}.  On the other hand, \Ito\/ differential equations are convenient for analysis as well as simulation. Thus, in formulating the Riemannian Langevin equation, it is necessary to first formulate the appropriate Stratonovich equation and then compute the deterministic \Ito\/--Stratonovich correction. A central observation in our work is that this correction term is due to curvature and is explicitly computable for several Riemannian geometries relevant to optimization~\cite{HIM22,IM,MY1,MY2}. 

\subsection{Riemannian Langevin Monte Carlo sampling schemes}
For the two metrics considered in this paper, the \Ito--Stratonovich correction due to curvature may also be computed explicitly, yielding the SDEs in Section~\ref{sec:SDE}. The rigorous analysis of these SDEs is presented in the companion paper~\cite{YZLMZ-theory}, and we focus on numerical algorithms in this paper. The Euler-Maruyama approximation to these SDEs yields the numerical sampling schemes listed below.

The SDEs also admit other numerical approximations. We have chosen the Euler-Maruyama schemes % that sample the Gibbs distribution 
because these schemes are fully explicit, simple to state, implement and numerically validate. They are generalizations of the popular unadjusted Langevin Monte Carlo for sampling in Euclidean spaces. Further, these schemes reduce to deterministic Riemannian gradient descent methods in the limit $\beta \to \infty$.

%Though there are many studies on Riemannian Langevin equations, to the best of our knowledge, none of existing work give or imply the same results in Section \ref{sec:SDE}.

%In particular, for the embedded geometry, our task reduces to the computation of the mean curvature of the embedding $\Snp \to \R^{n\times n}$. This follows from a general fact: when considering the extrinsic construction of Brownian motion on a manifold, the \Ito\/-Stratonovich correction is minus a half the mean curvature~\cite{IM}. 
%Another example of such curvature corrections relevant to stochastic optimization may be found in~\cite{MY1}.

%In summary, the main theoretical task is a formulation of the Riemannian Langevin equation, in both its \Ito\/ and Stratonovich variants. 

\subsubsection{Scheme E for the embedded geometry}
\label{sss:embed}
For the embedded manifold $(\mathcal{S}^{n,p}_+, g_E)$, the scheme is  
\begin{equation} \label{Euler_Maruyama_scheme_embedded}
 \resizebox{.99\hsize}{!}{$   X_{k+1} = \mbox{P}_{\mathcal{S}^{n,p}_+}\left[ X_k -\Delta t \,\grad \mathcal E(X_k) +  Q_k \left( \sqrt{\dfrac{2 \Delta t }{\beta}} \bmat{B_{11} & B_{12} \\ B_{12}^T & 0 } + \dfrac{\Delta t}{\beta}\sum\limits_{i=1}^p \frac{1}{\lambda_i} \bmat{0 & 0 \\ 0 &     I_{n-p}}  \right) Q_k^T \right]$},
\end{equation}
where $\mbox{P}_{\mathcal{S}^{n,p}_+}$ is the Euclidean projection to ${\mathcal{S}^{n,p}_+}$, and $X_k=Q_k\Lambda Q_k^T$ is the full SVD of $X_k =\mathcal{S}^{n,p}_+$ with eigenvalues $\lambda_1\geq \cdots\geq \lambda_p>0$.
The entries of $B_{12}$ are i.i.d.{} drawn from $\sqrt{\frac{1}{2}} \mathcal{N}(0,1)$. The entries of the symmetric $B_{11}$ are defined as follows: the diagonal entries are i.i.d. drawn from $ \mathcal{N}(0,1)$, and off-diagonal entries are  $b_{ij} = b_{ji} \sim \sqrt{\frac{1}{2}} \mathcal{N}(0,1)$. When $\beta=\infty$, equation~\eqref{Euler_Maruyama_scheme_embedded} reduces to $ X_{k+1} = \mbox{P}_{\mathcal{S}^{n,p}_+}\left( X_k -\Delta t\, \grad \mathcal E(X_k)\right)$, which is the Riemannian gradient descent on 
$(\mathcal{S}^{n,p}_+, g_E)$, see \cite{absil2008optimization, zheng2022riemannian}. We refer to \eqref{Euler_Maruyama_scheme_embedded} as Scheme E.

In this scheme, the  term $\dfrac{\Delta t}{\beta}\sum\limits_{i=1}^p \frac{1}{\lambda_i} \bmat{0 & 0 \\ 0 &     I_{n-p}}  $ in equation~\eqref{Euler_Maruyama_scheme_embedded} is the correction due to the mean curvature of the embedding of $\mathcal{S}^{n,p}_+\xhookrightarrow{}  \mathbbm R^{n\times n}$. 
\subsubsection{Scheme BW for the Bures-Wasserstein metric}
\label{sss:bw}
For the quotient manifold $( \mathbb{R}^{n\times p}_*/\mathcal{O}_p, g_{BW})$, the scheme is   
\begin{equation}\label{Euler_Maruyama_scheme_quotient}
    Y_{k+1} = Y_{k} - \Delta t   2 \nabla \mathcal{E} (Y_kY_k^T) Y_k + \sqrt{\frac{2\Delta t}{\beta}} B_k + \frac{\Delta t}{\beta} U_k \bmat{\sum\limits_{j: j\neq i}\frac{\sigma_i}{\sigma_i^2 + \sigma_j^2}}_{ii} V_k^T, 
\end{equation}
where $B_k$ is $n$-by-$p$ matrix with entries being i.i.d.{} standard Gaussian, 
$Y_k = U_k \Sigma_k V_k^T \in \mathbb R^{n\times p}$ is the compact SVD with singular values $\sigma_i$, and $\bmat{\sum\limits_{j: j\neq i}\frac{\sigma_i}{\sigma_i^2 + \sigma_j^2}}_{ii}$ is the diagonal matrix whose $i$-th diagonal entry is $\sum\limits_{j: j\neq i}\frac{\sigma_i}{\sigma_i^2 + \sigma_j^2}$. We refer to \eqref{Euler_Maruyama_scheme_quotient} as Scheme BW.
The Riemannian Langevin Monte Carlo scheme~\eqref{Euler_Maruyama_scheme_quotient} can be viewed as  a natural extension of Burer-Monteiro gradient descent method 
\begin{equation}
    Y_{k+1} = Y_{k} - \Delta t\,   2 \nabla \mathcal{E} (Y_kY_k^T) Y_k ,
    \label{BM-GD}
\end{equation}
which is the simplest low-rank gradient descent method for minimizing $\mathcal E(X)$ under the constraint $X\in \mathcal S_+^{n,p}$. It is clear that as $\beta\to\infty$, \eqref{Euler_Maruyama_scheme_quotient} reduces to \eqref{BM-GD}.
The Burer-Monteiro gradient descent method is equivalent to a Riemannian gradient descent method on the quotient manifold $\mathbb{R}^{n\times p}_*/\mathcal{O}_p$ with Bures-Wasserstein metric, see \cite{zheng2022riemannian}.

\subsubsection{Gibbs distribution sampling and numerical validation}
%The reader should recognize that while both schemes are for matrices on the same {\em differentiable\/} manifold $\Snp$, the two schemes correspond to two different {\em Riemannian\/} manifolds. In particular, 
While the Gibbs distribution always has the same density function $e^{-\beta\mathcal E}$ with respect to $\rhoref$, the reference density $\rhoref$ depends on the metric. Thus, 
the two schemes \eqref{Euler_Maruyama_scheme_embedded} and \eqref{Euler_Maruyama_scheme_quotient}, generate samples for two different probability distributions. In order to validate our schemes, we choose energy functions that allow an explicit computation of these densities for both metric.  These energy functions yield matrix integrals of independent analytic interest. They also allow side-to-side benchmarking for different Gibbs samplers on $\Snp$. We demonstrate the efficiency of sampling from these Gibbs distributions numerically. Further analysis on convergence to equilibrium as $t\to \infty$ using the Bakry-Emery criterion is considered in the companion paper~\cite{YZLMZ-theory}. 

Finally, while we do not discuss their convergence and efficiency approximating SDE as the step-size $\Delta t \to 0$; this is possible following existing approximation results \cite{MR1214374, NEURIPS2022_27c852e9, li2023riemannian}.

%While there are several examples of Another focus of this paper is to demonstrate the numerical convergence of  simple schemes \eqref{Euler_Maruyama_scheme_embedded} and \eqref{Euler_Maruyama_scheme_quotient}. One would need at least the probability density distribution of some random variable  on  the manifold, which is possible for special manifolds such as spheres, torus and orthogonal group \cite{diaconis2013sampling, zappa2018monte, ge2020efficient}.  For$\mathcal S_+^{n,p}$, it is  in general much more difficult to compute  the Gibbs distribution due to the Riemannian volume form  $\ud V$. In Section \ref{sec:example}, we construct a few examples on  $\mathcal S_+^{n,p}$ with analytical formulae for both metrics, which is another contrubtion. We   test the two schemes numerically in Section \ref{sec:tests}. 

\subsection{Some applications and related work}
\label{subsec:applications}

% [GM: this section should reflect representative applications and cite leading books and papers in the areas. We could also just keep it short since Gibbs sampling, SDP, ML are already very big areas.] \jl{agreed, needs to revisit later}

\subsubsection{Applications of PSD matrices}
Positive semi-definite (PSD) fixed rank matrices arise in many problems such as distance matrices \cite{tasissa2018exact} and covariance matrices in statistics, and
have been   used in   applications  including  kernels in machine learning \cite{meyer2011regression}, semidefinite optimization \cite{burer2005local}, quantum information, etc.
  Riemannian optimization algorithms over $\mathcal{S}^{n,p}_+$ under different metrics have been well studied, e.g., see \cite{vandereycken2009embedded,huang2017solving, massart2019curvature, zheng2022riemannian} and references therein. 

\subsubsection{Langevin dynamics and Monte Carlo schemes on manifolds}\ 

There is an extensive literature on Langevin dynamics in  statistics and related areas, with interest in
nonconvex optimization \cite{cheng2018underdamped, cheng2020stochastic}, as well as machine learning such as  generative models \cite{du2019implicit}. 
 
In recent years, there has been interest in studying Langevin diffusion and Monte Carlo Markov Chain (MCMC) schemes on manifolds  \cite{ciccotti2005blue, ciccotti2008projection, girolami2011riemann, pmlr-v22-brubaker12, byrne2013geodesic, zappa2018monte, moitra2020fast, ge2021efficient, leake2021sampling, li2023riemannian}.  
In this paper, we are interested in Riemannian Langevin Monte Carlo schemes on $\mathcal{S}^{n,p}_+$.

In the statistics literature, manifold Langevin schemes have been studied in~\cite{girolami2011riemann,byrne2013geodesic}. However, these schemes apply  only to simpler embedded manifolds $\mathcal M\subset\mathbbm R^{n}$ with explicit geodesics such as the sphere and Stiefel manifolds. 
The above schemes do not directly apply to the manifold $\mathcal S_+^{n,p}$, even for the embedded geometry.
In \cite{zappa2018monte}, a sampling scheme using projection to surface is constructed; however,  this is not a Langevin scheme.

%To construct a simple Riemannian Langevin scheme,
%the first difficulty one encounters is the Brownian $\bd B_t^{g,\beta}$ on a complicated manifold  $(\mathcal S_+^{n,p},g)$.
%In the literature of probability theory and stochastic processes,
%SDEs on a Riemannian manifold have been well studied  \cite{stroock2000introduction, MR1882015},
%in which the Brownian motion $\bd B_t^{g,\beta}$ is usually represented in Stratonovich stochastic integral. The Stratonovich stochastic integral is an implicit mid-point rule, which is not very useful if we prefer to construct simple explicit Langevin schemes.
%It is possible to rewrite $\bd B_t^{g,\beta}$ equivalently as an It\^o Brownian motion term with a deterministic geometric correction [{\it cite Govind's paper}]. In a separate paper [{\it cite analysis paper}], we give the full proof of these SDEs on $\mathcal{S}^{n,p}_+$. 

%When $\bd B_t^{g,\beta}$ is given in It\^o form, we can use explicit discretizations, e.g., the simple Euler-Maruya method for SDE to obtain simple Riemannian Langevin schemes. 

In general, a Langevin scheme can be used for either optimization \cite{xu2020global, li2023riemannian}, or Monte Carlo type numerical integration, which is common in Bayesian statistic. 
For optimization, stochastic optimization by Langevin dynamics with simulated annealing is an established approach    \cite{liu2001monte}. In \cite{cheng2018underdamped}, underdamped Langevin schemes are shown to be much more efficient than the overdamped case \eqref{scheme-LD}.
For sampling, Metropolis-adjusted Langevin algorithm \cite{girolami2011riemann} is often used. 
For simplicity, we focus on the simple schemes \eqref{Euler_Maruyama_scheme_embedded} and \eqref{Euler_Maruyama_scheme_quotient} without considering any of simulated annealing, underdamped Langevin, or Metropolis-adjustment, to which it is possible to extend our schemes. 
Though the Riemannian optimization on $\mathcal S_+^{n,p}$ can be easily extended to Hermitian PSD matrices of fixed rank \cite{zheng2022riemannian}, we remark that such an extension for Langevin dynamics would be significantly different. 

\subsection{Organization of the paper}
In Section~\ref{sec:SDE}, we state the explicit formulae for the SDE \eqref{sde1} and Gibbs measure on the manifold $\mathcal S_+^{n,p}$ under two metrics $g_E$ and $g_{BW}$. We then derive the  schemes \eqref{Euler_Maruyama_scheme_embedded} and \eqref{Euler_Maruyama_scheme_quotient} in Section \ref{sec:scheme}. The energy functions and Gibbs distributions used to benchmark the schemes are presented in Section~\ref{sec:example}. The numerical results are studied in Section~\ref{sec:tests}.
%  Concluding remarks are given Section \ref{sec:remark}.

\section{Riemannian Langevin equations on $\mathcal{S}^{n,p}_+$}
\label{sec:SDE}
In this section, we state the \Ito\/ form of the Riemannian Langevin equation~\eqref{sde1} for both Riemannian geometries studied in this paper. The theoretical basis for these SDEs is discussed at greater depth in~\cite{YZLMZ-theory}. The main ideas are as follows: (a) the abstract theory of Brownian motion on Riemannian manifolds is used to define the Riemannian Langevin equation in Stratonovich form for the metrics $g_E$ and $g_{BW}$ on $\Snp$; (b) the \Ito\/-Stratonovich conversion rule is used to compute the associated \Ito\/ form of these SDEs and it is observed that the \Ito-Stratonovich\/ correction term corresponds to mean curvature.  
This approach yields the SDEs below. These SDEs are used to develop numerical schemes in Section~\ref{sec:scheme}.

%For the embedded geometry, we use the observation that the correction term is minus one half mean curvature~\cite{IM}, reducing our task to the computation of the mean curvature of the embedding. For the quotient geometry, we use the observation (as in~\cite{HIM22,MY2}) that the correction is the projection of the 

\subsection{The Riemannian Langevin equation for 
embedded geometry $(\mathcal S_+^{n,p}, g_{E})$}
Let $X\in S_{+}^{n,p}$ have the compact SVD $X=U\Lambda U^T$ with $U\in \mathbbm R^{n\times p}$. Let $U_\perp\in \mathbbm R^{n\times (n-p)}$ be a matrix with columns orthonormal to columns of $U$. The tangent space of $S_{+}^{n,p}$ at $X=U\Lambda U^T\in S_{+}^{n,p}$ is given by \cite{vandereycken2009embedded, zheng2022riemannian}:
\begin{equation}
    \label{tangent-space}
    T_XS_{+}^{n,p}=\left\{\begin{bmatrix}
        U & U_\perp
    \end{bmatrix}\begin{bmatrix}
        H & K^T\\
        K & 0
    \end{bmatrix}\begin{bmatrix}
        U^T \\ U^T_\perp
    \end{bmatrix}: \forall K\in \mathbbm R^{(n-p)\times p},\forall H\in \mathbbm R^{p\times p}, H^T=H  \right\}. 
\end{equation}
The induced metric $g_E$ by the embedding $\mathcal{S}^{n,p}_+\xhookrightarrow{}  \mathbbm R^{n\times n}$ is then defined as $$g_E(A,B)=\Tr(A^TB),\quad \forall A, B\in T_XS_{+}^{n,p},$$ which is the Frobenius inner product for two matrices.

%Notice that \eqref{sde-stratonovich} involves  orthonormal basis of $\mathfrak X(T\mathcal M)$ and Stratonovich stochastical differential, for both of which it would be difficult to construct simple numerical schemes. 

Equation~\eqref{sde1} describes the evolution of a point $\bd X_t \in \Snp$ in abstract terms. We now rewrite it in a simpler equivalent form describing the evolution of the entries of the matrix entries $\{(X_t)_{ij}\}_{i,j=1}^n$ representing $\bd X_t$. Let us write $X=U\Lambda U^T$ for the compact singular value decomposition (SVD) of $X$. We further assume that the singular values $\Lambda=\mathrm{diag}(\lambda_1,...,\lambda_p)$ are written in decreasing order. We suppress the subscript $t$ in the following equations, though the reader should note that $U$ and $\Lambda$ depend on $X_t$. 

Then we find that the law of $X_t$ is determined by the \Ito\/ differential equation
\begin{align}\label{sde2}
        \ud X_t=-\grad \mathcal E(X_t)\ud t+\sqrt{\frac2\beta} \ud W^{n,p, X_t}_t+\frac1\beta H(X_t)\ud t.
\end{align}
In this equation, the stochastic forcing $W^{n,p, X_t}_t$ is the orthogonal projection of white noise in $\R^{n\times n}$ onto $T_{X_t}\Snp$. Precisely, given $W^i_t $ for $1\leq i\leq n$ and $W^{i,j}_t$ for $1\leq i< j \leq n$  independent standard one-dimensional Wiener process, we set 
    \begin{align}
    \label{eq:embed-sde-noise}
        \resizebox{.99\hsize}{!}{$ \ud W^{n,p, X_t}_t=\begin{bmatrix}
            U & U_\perp 
        \end{bmatrix} \begin{bmatrix}
            \ud W^{1}_t&\cdots&\frac1{\sqrt2}\ud W^{1,p}_t&\frac1{\sqrt2}\ud W^{1,p+1}_t&\cdots&\frac1{\sqrt2}\ud W^{1,n}_t\\
            \vdots&\ddots&\vdots&\vdots&\ddots&\vdots\\
            \frac1{\sqrt2}\ud W^{1,p}_t&\cdots&\ud W^{p}_t&\frac1{\sqrt2}\ud W^{p,p+1}_t&\cdots&\frac1{\sqrt2}\ud W^{p,n}_t\\
            \frac1{\sqrt2}\ud W^{1,p+1}_t&\cdots&\frac1{\sqrt2}\ud W^{p,p+1}_t&0&\cdots&0\\
            \vdots&\ddots&\vdots&\vdots&\ddots&\vdots\\
            \frac1{\sqrt2}\ud W^{1,n}_t&\cdots&\frac1{\sqrt2}\ud W^{p,n}_t&0&\cdots &0
        \end{bmatrix}\begin{bmatrix}
            U^T\\  U_\perp^T 
        \end{bmatrix}\nonumber, $}
     \end{align}   
The term $H(X_t)$ is the mean curvature of the embedding $\Snp \to \R^{n\times n}$. We adopt the convention in geometric analysis: the mean curvature is defined as the trace of the second fundamental form of the embedding. Explicitly, we have
          \begin{align}
          \label{eq:embed-sde-mc}
        H(X_t)=\left(\sum\limits_{i=1}^p \frac1{\lambda_{i}}\right)\begin{bmatrix}
            U & U_\perp 
        \end{bmatrix}\begin{bmatrix}
            0_{p\times p}&0_{p\times(n-p)}\\
            0_{(n-p)\times p}& I_{n-p}
        \end{bmatrix}\begin{bmatrix}
            U^T\\  U_\perp^T 
        \end{bmatrix}.
    \end{align}
The following feature of equation~\eqref{sde2} is fundamental. The stochastic forcing is the naive projection of white noise in the ambient space $\mathbb{R}^{n\times n}$ onto $T_{X_t}\Snp$. Intuitively, when one uses the Euler-Maruyama discretization, the role of this term is to update $X_t$ by taking unbiased random steps in any direction in the tangent space. However, \Ito\/ calculus has a subtle interplay with the geometry of the embedding, and in order to keep $X_t$ on the manifold $\Snp$, it is necessary to include the correction term given by the mean curvature. 

%Several variants of this insight have been considered by two of the authors in recent work~\cite{HIM22,IM,MY1,MY2}.

\subsection{The Riemannian Langevin equation for 
Bures-Wasserstein geometry $(\mathcal S_+^{n,p}, g_{BW})$}
 The manifold ${S}^{n,p}_+$ can also be viewed as a quotient manifold $\mathbb{R}^{n\times p}_*/\mathcal{O}_p$, for which the noncompact Stiefel manifold $\mathbb{R}^{n\times p}_*$ is called the \textit{total space}.
Denote the natural projection as
$$\pi: \mathbb{R}^{n\times p}_* \rightarrow \mathbb{R}^{n\times p}_*/\mathcal{O}_p.$$ 
For any $Y\in \mathbb R_*^{n\times p}$,    the equivalence class containing $Y$ is
$$[Y]=\pi^{-1}(\pi(Y)) =\left\{ YO \mid O\in \mathcal{O}_p \right\},$$
   which is an embedded submanifold of  $\mathbb{R}^{n\times p}_*$ (see e.g., \cite[Prop.~3.4.4]{absil2008optimization}). The tangent space of $[Y] $ at $Y$ is therefore a subspace of $T_Y \mathbb{R}^{n\times p}_*$ called the \textit{vertical space} at $Y$, and is denoted by $\mathcal{V}_Y=\left\{Y\Omega\mid \Omega^T  = - \Omega, \Omega \in \mathbb{R}^{p\times p} \right\}$, see \cite{zheng2022riemannian}.  

Define $$\begin{aligned} \theta : \mathbb{R}^{n\times p}_* & \rightarrow \mathcal{S}^{n, p }_+\\
 Y &\mapsto YY^T.\end{aligned}$$  Then $\theta$ is invariant under the equivalence relation  and induces a bijection $\Tilde{\theta}$ on $\mathbb{R}^{n\times p}_*/\mathcal{O}_p$  such that $\theta = \Tilde{\theta} \circ \pi$. 
For any  function $\mathcal E(X)$ defined on $\mathcal{S}^{n,p}_+$, there is a  function $F$ defined on $\mathbb{R}^{n\times p}_*$ that induces $\mathcal E$: for any $X = YY^T \in \mathcal{S}^{n,p}_+$, $F(Y) := \mathcal E\circ \theta (Y) = \mathcal E(YY^T)$. This is summarized in the diagram below:
\[
  \begin{tikzcd}[baseline=\the\dimexpr\fontdimen22\textfont2\relax]
     \mathbb{R}^{n\times p}_* \arrow[rd,"\theta:= \tilde{\theta}\circ \pi",dashrightarrow] \arrow[d,"\pi"] & & \\
     \mathbb{R}^{n\times p}_*/\mathcal{O}_p \arrow[r,leftrightarrow,"\tilde{\theta}"]  &\mathcal{S}^{n,p}_+ \arrow[r,"\mathcal E"] & \mathbb{R}
  \end{tikzcd}
\] 
In particular, $\mathcal S_+^{n,p}$ is diffeomorphic to
$\mathbb{R}^{n\times p}_*/\mathcal{O}_p $ under $\tilde \theta$, see
 \cite{zheng2022riemannian}.  
 For any $Y\in \mathbb{R}_*^{n\times p}$, the flat metric for the total space $\mathbb{R}^{n\times p}_*$, correction term corresponds to mean curvature.  $$g(a,b)=\Tr(a^T b),\forall a,b\in T_Y\mathbb{R}_*^{n\times p}=\mathbb{R}^{n\times p}$$ induces a metric on the quotient manifold $\mathbb{R}^{n\times p}_*/\mathcal{O}_p $, which is called Bures-Wasserstein metric, see \cite{massart2019curvature,massart2020quotient, zheng2022riemannian}.
Another way to understand the Bures-Wasserstein metric at $X\in \mathcal{S}^{n,p}_+\eqsim \mathbb{R}^{n\times p}_*/\mathcal{O}_p$ is 
  via the map $\theta$: 
  % \jl{is it the convention that we suppress the base of the differential $d \beta$ in notation?} \ty{Yes. Sometimes one also uses $\beta_*$ instead of $\ud\beta$ to denote the pushforward map. I changed $\beta$ to $\theta$ as $\beta$ is also the notation for inverse temperature.}
  \begin{align}
    &g_{BW}(A,B)=\Tr(ab^T)&&\forall A,B\in T_X\mathcal{S}^{n,p}_+, a,b\in T_Y\mathbb R_*^{n\times p}\nonumber\\
    &&&\st \ud \theta(Y)[a]=A, \ud \theta(Y) [b]=B, a,b\in \mathrm{Ker}(\ud \theta(Y) )^\perp
\end{align}
where $X$ has decomposition $X=YY^T$, $\ud \theta(Y)[a]=Ya^T+aY^T$ is the differential of $\theta$ at $Y$, and $a\in\mathrm{Ker}(\ud \theta (Y))^\perp\Leftrightarrow Y^Ta=a^TY$.

%\subsection{The Riemannian Langevin equation under the Bures-Wasserstein metric}
The Riemannian Langevin equation is now determined by the geometry of Riemannian submersion. We must obtain an \Ito\/ differential equation for $Y_t$, such that $X_t=Y_tY_t^T$ is a matrix that has the same law as the solution to \eqref{sde1} in $(\mathcal{S}^{n,p}_+,g_{BW})$. 

In comparison with equation~\eqref{sde2}, we see that the natural choice for white noise driving $Y_t$ is white noise in $\mathbb{R}^{n\times p}$. This is the stochastic differential $\ud W_t$, where $W_t = \{W_t^{ij}\}_{1\leq i \leq n, 1 \leq j \leq p}$ consists of $np$ independent standard one-dimensional Wiener processes. However, as in equation~\eqref{sde2} we must include a deterministic correction. This correction corresponds to mean curvature again, but in a more subtle way than~\eqref{sde2}. The equivalence class of $Y$ such that $X=YY^T$ is a group orbit of $\mathcal{O}_p$ embedded within $\mathbb{R}^{n\times p}$. The logarithm of the volume of this group orbit constitutes a natural Boltzmann entropy. It may be computed explicitly, and we find
\begin{equation}
\label{eq:entropy}
S(Y)=\frac1{2}\sum\limits_{i=1}^p\sum\limits_{j=i+1}^p \log(\sigma_i^2+\sigma_j^2)
\end{equation}
where $\{\sigma_i\}_{i=1}^p$ are  singular values of $Y$. It is known that  $\nabla S(Y)$ is the mean curvature of the group orbit in $\mathbb{R}^{n\times p}$~\cite[p.3505]{Pacini}.

We then have the following \Ito\/ differential equation for $Y_t$ such that $X_t=Y_tY_t^T$ has the same law as the solution to~\eqref{sde1}.
    \begin{align}
    \label{SDE-BW}
        \ud Y_{ij}=&-\frac{\partial \mathcal E(YY^T)}{\partial Y_{ij}}\ud t+\sqrt{\frac2\beta}\ud W^{ij}_t-\frac{1}{\beta}\frac{\partial  S(Y)}{\partial Y_{ij}}\ud t, &&1\leq i\leq n,1\leq j\leq p.
    \end{align}
The correction term can be explicitly computed using the following
\begin{lem}
If $Y\in\mathbb\mathbb R^{n\times p}_*$ has SVD as $Y=Q\Sigma P^T$
with singular values $\sigma_i$, then 
    the gradient of the correction term $ S$ is given by
    $\nabla S(Y)=Q \tilde \Sigma P^T$ where $\tilde \Sigma$ is a diagonal matrix with diagonal entries 
    $\sum_{j\neq 1}\frac{\sigma_1}{\sigma_1^2+\sigma_j^2}, \sum_{j\neq 2}\frac{\sigma_2}{\sigma_2^2+\sigma_j^2},\cdots,  \sum_{j\neq p}\frac{\sigma_p}{\sigma_p^2+\sigma_j^2}.$
\end{lem}

\section{Two Riemannian Langevin Monte Carlo   schemes}
\label{sec:scheme}
% \jl{we'd better give a uniform name for these, perhaps Reimannian Langevin Monte Carlo schemes?}
To get a simple Riemannian Langevin Monte Carlo sampling scheme, we only consider convenient discretization and approximation methods, which can be easily and efficiently implemented. 
For the Brownian motion term, we consider the most straightforward and simplest
discretization of the SDEs \eqref{sde2} and \eqref{SDE-BW}, i.e., the
Euler-Maruyama type discretization.

One extra complication from the manifold constraint is how to approximate the exponential map. 
For optimization algorithms on Riemannian manifolds \cite{absil2008optimization}, retraction, which is at least a first order
approximation to the exponential map, is often used. For instance, for approximating an ODE $\frac{\ud}{\ud t}\bd X=-\grad \mathcal E(\bd X)$ on a manifold $\mathcal M$, with any retraction operator $\mathcal R_{\mathcal M}$ mapping to $\mathcal M$, a simple forward Euler type approximation, or equivalently the Riemannian gradient descent method, is given by
\[\bd X_{k+1}=\mathcal R_{\mathcal M}[\bd X_{k+1}-\Delta t\,\grad \mathcal E(\bd X_k)].\]

In particular, when combining the  Euler-Maruyama type discretization for SDE and the simple Riemannian gradient descent by retraction, we get the two simple Riemmanian Langevin Monte Carlo schemes as follows.

\subsection{Scheme E for the embedded geometry}

\subsubsection{The  Riemannian gradient} 
For a given energy function $\mathcal E(X)$,
its Riemannian gradient $\grad \mathcal E(X)$ of  at $X\in \mathcal S_+^{n,p}$, is the Euclidean projection of the Euclidean gradient $\nabla \mathcal E(X)\in\mathbb R^{n\times n}$ defined as $[\nabla \mathcal E(X)]_{ij}=\frac{\partial}{\partial X_{i}} \mathcal E(X)$, onto the tangent space $T_X \mathcal S_+^{n,p}$, see \cite{absil2008optimization,vandereycken2009embedded,zheng2022riemannian}. 
It is straightforward to verify that $\nabla \mathcal E(X)$ is a symmetric matrix for any differentiable $\mathcal E$ and any $X\in \mathcal S_+^{n,p}$.
For any given $X\in \mathcal S_+^{n,p}$, let $X=U\Lambda U^T$ be its compact SVD. 
Let $P_U=UU^T$ and $P_{U_\perp}=U_\perp U_\perp^T=I-UU^T$.
By derivations in \cite{zheng2022riemannian}, $\grad \mathcal E(X)$ can be computed and represented as 
 \begin{align*}
      \grad \mathcal E(X)&=\begin{bmatrix}
        U & U_\perp
    \end{bmatrix}\begin{bmatrix}
        U^T\nabla \mathcal E(X)U & U^T\nabla \mathcal E(X)U_\perp\\
        U_\perp^T\nabla \mathcal E(X)U & 0
    \end{bmatrix}\begin{bmatrix}
        U^T \\ U^T_\perp
    \end{bmatrix}\\
    &=P_U\nabla \mathcal E(X)P_U+P_{U_\perp}\nabla \mathcal E(X)P_U+P_U\nabla \mathcal E(X)P_{U_\perp}.
 \end{align*} 
 The compact implementation of computing $\grad \mathcal E(X)$  is given in Algorithm \ref{alg:grad_embedded}.

\begin{algorithm}[H]
\caption{Compact computation of the Riemannian gradient $\grad \mathcal E(X)$}\label{alg:grad_embedded}
\begin{algorithmic}[1]
\Require The compact SVD of $X\in {S}^{n,p}_+$: $X = U\Lambda U^T$
\Ensure $\grad \mathcal E(X) = UHU^T + U_p U^T + U U_p^T \in T_X{S}^{n,p}_+$ 
\par $T \leftarrow \nabla \mathcal E(X)U$ 
\par $H \leftarrow U^TT$ 
\par $U_p \leftarrow T - UH$  
\end{algorithmic}
\end{algorithm}

\subsubsection{The
retraction by projection}
Let $\mathcal{S}^{n\times n}$ denote symmetric matrices, then the Euclidean projection $\mbox{P}_{{S}^{n,p}_+}: \mathcal{S}^{n\times n}\longrightarrow \mathcal S_+^{n,p}$
 is  a convenient retraction operator, see \cite{absil2008optimization,vandereycken2009embedded,zheng2022riemannian}. 
A straightforward implementation is given in Aglorithm \ref{alg:retraction_embedded}.

\begin{algorithm}[H]
\caption{Computation of the retraction $P_{{S}^{n,p}_+}(X+Z)$}\label{alg:retraction_embedded}
\begin{algorithmic}[1]
\Require the compact SVD of $X$: $X = U\Lambda U^T \in {S}^{n,p}_+$,  $Z \in \mathcal{S}^{n\times n}$. 
\Ensure $P_{{S}^{n,p}_+}(X+Z) = Q_+\Lambda_+ Q_+^T\in {S}^{n,p}_+$. 
\par $(Q_+,\Lambda_+) = \mbox{svd}(X+Z)$ 
\par $U_+ \leftarrow Q_+(:,1:p)$ \quad $\Lambda_+ \leftarrow \Lambda_+(1:p,1:p)$
\end{algorithmic}
\end{algorithm}

\subsubsection{A Riemannian Langevin Monte Carlo scheme}
For approximating the SDE \eqref{sde2} on $(S_{+}^{n,p}, g_E)$,
with the retraction operator and Euler-Maruyama method for SDE, 
we have the   scheme \eqref{Euler_Maruyama_scheme_embedded}, which can be 
more explicitly written as
\begin{equation}\label{Euler_Maruyama_scheme_embedded2}
   \resizebox{.99\hsize}{!}{$   X_{k+1} = \mbox{P}_{{S}^{n,p}_+}\left(\begin{bmatrix}
        U & U_\perp
    \end{bmatrix}\begin{bmatrix}
        \Lambda-\Delta t U^T\nabla \mathcal E(X_k)U+\sqrt{\frac{2 \Delta t }{\beta}}B_{11} & -\Delta t U^T\nabla \mathcal E(X_k)U_\perp+\sqrt{\frac{2 \Delta t }{\beta}}B_{12}\\
        -\Delta t U_\perp^T\nabla \mathcal E(X_k)U+\sqrt{\frac{2 \Delta t }{\beta}}B_{12}^T & \frac{\Delta t}{\beta} \sum\limits_{i=1}^p \frac{1}{\lambda_i} I_{n-p}
    \end{bmatrix}\begin{bmatrix}
        U^T \\ U^T_\perp
    \end{bmatrix} \right),$}
\end{equation}
where 
$X_k = U \Lambda U^T$ is the compact SVD of $X_k\in {S}^{n,p}_+$ with eigenvalues $\lambda_1\geq \lambda_2\geq\cdots\geq \lambda_p>0$.
The third term in the right hand side is the white noise term in the tangent space $T_{X_k} \mathcal S_+^{n,p}$.
Entries of $B_{12}\in\mathbbm R^{p\times(n-p)}$ are i.i.d drawn from $\sqrt{\frac{1}{2}} \mathcal{N}(0,1)$, and $B_{11}\in\mathbbm R^{p\times p}$ are defined as follows. 
\begin{equation}
    B_{11} = \begin{bmatrix}
 \mathcal{N}(0,1) &         &         &  \\
       & \ddots  &  b_{ij}  &    \\
       &  b_{ji} & \ddots &    \\
       &         &        &  \mathcal{N}(0,1)   \\
\end{bmatrix}
\end{equation}
with $b_{ij} = b_{ji} \sim \sqrt{\frac{1}{2}} \mathcal{N}(0,1)$. 
The implementation details of the scheme \eqref{Euler_Maruyama_scheme_embedded} are given as follows in the Algorithm \ref{alg:Embedded_Langevin_Dynamics}.

\begin{algorithm}[H]
\caption{The Riemannian Langevin Monte Carlo scheme \eqref{Euler_Maruyama_scheme_embedded} for $(\mathcal S_+^{n,p}, g_E)$}\label{alg:Embedded_Langevin_Dynamics}
\begin{algorithmic}[1]
\Require initial iterate $X_1 \in \mathcal{S}^{n,p}_+$; full SVD of $X_1$: $X_1 = Q_1 \Lambda_1 Q_1^T$ 
\For{$k =1,2,\dots, N$}
    \State{Compute Riemannian gradient}
    \par\hskip\algorithmicindent $\xi_k := \grad \mathcal E(X_k)$ 
    \Comment{ See Algorithm \ref{alg:grad_embedded}} 
    \State{Compute noise term}
    \par\hskip\algorithmicindent $B = \sqrt{\frac{2 \Delta t }{\beta}} \bmat{B_{11} & B_{12} \\ B_{12}^T & 0 } + \frac{\Delta t}{\beta} \sum\limits_{i=1}^p \frac{1}{\lambda_i} \bmat{0 & 0 \\ 0 &    I_{n-p}} $
    \State Obtain the new iterate by retraction
    \par\hskip\algorithmicindent $X_{k+1} = \mbox{P}_{\mathcal{S}^{n,p}_+}(X_k-\Delta t \xi_k + Q_k B Q_k^T)$ \Comment{ See Algorithm \ref{alg:retraction_embedded}} 
\EndFor
\end{algorithmic}
\end{algorithm}

 \begin{rem}
The  mean curvature correction term is necessary for avoiding rank deficient samples in the following sense. 
A sampling scheme on $\mathcal{S}^{n,p}_+$ might generate a sample $X$ with a rank numerically close to $p-1$, and the mean curvature correction term in the scheme \eqref{Euler_Maruyama_scheme_embedded} would be huge if $\lambda_p\to 0$, thus it will force iterate $X_k$ to stay away from the boundary of $\mathcal{S}^{n,p}_+$.  
\end{rem}
\begin{rem}
Notice that the complexity of computing SVD of $X+Z$ in Algorithm \ref{alg:retraction_embedded} would be $\mathcal O(n^3)$ in a naive implementation. For a Riemannian gradient method, if $Z\in T_{X_k}\mathcal{S}^{n,p}_+$, a compact implementation of computing $P_{\mathcal{S}^{n,p}_+}(X+Z)$ in \cite{zheng2022riemannian} is only $\mathcal O(np^2)+\mathcal O(p^3)$,
which is no longer possible for the Langevin Monte Carlo scheme   \eqref{Euler_Maruyama_scheme_embedded} due to the mean curvature correction term in the normal space.
On the other hand, if Lanczos type algorithm is used for computing to top $p$ eigen-componenes of $X+Z$, it seems possible to explore the special structure in \eqref{Euler_Maruyama_scheme_embedded2} to find a more efficient implementation, but we do not consider a more compact implementation in this paper. 
 \end{rem}

\subsection{Scheme BW for the Bures-Wasserstein metric}

\subsubsection{The Riemannian gradient  and a simple retraction operator}\ 

Given a smooth energy function $\mathcal E(X)$ defined  on $\mathcal{S}^{n,p}_+$, the corresponding function $h$ on $\mathbb{R}^{n\times p}_*/\mathcal{O}_p$ satisfies
\begin{equation}\label{eqn:cost_function_quotient}
\begin{aligned}
    h: \mathbb{R}^{n\times p}_*/\mathcal{O}_p & \rightarrow \mathbb{R}\\
     \pi(Y) & \mapsto \mathcal E(\tilde{\beta}(\pi(Y))) = \mathcal E(\beta(Y)) = \mathcal E(YY^T). 
\end{aligned}
\end{equation}
Observe that the function $F(Y):=\mathcal E(YY^T)$ satisfies  $F(Y)= h\circ \pi(Y) = \mathcal E \circ \beta(Y)$. 
The Riemannian gradient of $h$ at $\pi(Y)$ is a tangent vector in $T_{\pi(Y)}\mathbb{R}^{n\times p}_*/\mathcal{O}_p$ . The next theorem is given in \cite[Section 3.6.2]{absil2008optimization}, showing that the horizontal lift of $\grad h(\pi(Y))$ can be obtained from the Riemannian gradient of $F$ defined on $\mathbb{R}^{n\times p}_*$. 
\begin{thm}\label{thm:lift_gradient_quotient} The horizontal lift of the gradient of $h$ at $\pi(Y)$ is the Riemannian gradient of $F$ at $Y$. That is, 
\[
\overline{\grad h(\pi(Y)) }_Y  = \grad F(Y).
\]
\end{thm} 
For the Bures-Wasserstein metric, the following result is proven in \cite{zheng2022riemannian}: 
\begin{prop}\label{thm:gradF(Y)}
Let $\mathcal E$ be a smooth real-valued function defined on $\mathcal{S}^{n,p}_+$ and let $F: \mathbb{R}^{n\times p}_* \rightarrow \mathbb{R}: Y \mapsto \mathcal E(YY^T)$. Assume $YY^T = X$. Then the Riemannian gradient of $F$ is given by
\[
    \grad F(Y) = 
     2\nabla \mathcal E(YY^T)Y
\]
where $\nabla \mathcal E(\cdot)$ is the gradient of $\mathcal E$ w.r.t. $X$.
\end{prop}

In  \cite[Prop.~A.8]{massart2020quotient}, the relationship between the horizontal lifts of the  quotient tangent vector $\xi_{\pi(Y)}$ lifted at different representatives in $[Y]$ is given:
\begin{lem}\label{lem:horizontal_lift}
Let $\eta$ be a vector field on $\mathbb{R}^{n\times p}_*/\mathcal{O}_p$, and let $\bar{\eta}$ be the horizontal lift of $\eta$. Then for each $Y\in \mathbb{R}^{n\times p}_*$, we have $$\bar{\eta}_{YO} = \bar{\eta}_YO$$ for all $O\in \mathcal{O}_p$. 
\end{lem}

 The retraction on the quotient manifold $\mathbb{R}^{n\times p}_*/\mathcal{O}_p$ can be defined using the retraction on the total space $\mathbb{R}^{n\times p}_*$.  
For any $A \in T_Y \mathbb{R}^{n\times p}_*$ and a step size $\tau>0$,
\[
\overline{R}_Y(\tau A):= Y + \tau A, 
\]
is a retraction on $\mathbb{R}^{n\times p}_*$ if  $Y+\tau A$ remains full rank, which is ensured for small enough $\tau$. 
Then Lemma \ref{lem:horizontal_lift} indicates that $\overline{R}$ satisfies the conditions of \cite[Prop.~4.1.3]{absil2008optimization}, which implies that 
\begin{equation}\label{eqn:retraction_quotient}
    R_{\pi(Y)}(\tau \eta_{\pi(Y)}) := \pi(\overline{R}_Y(\tau  \overline{\eta}_Y)) = \pi(Y+\tau \overline{\eta}_Y)
\end{equation}
defines a retraction on the quotient manifold $\mathbb{R}^{n\times p}_*/\mathcal{O}_p$
for a small enough step size $\tau >0.$

Finally, we give an example of what these results imply by considering the Riemannian gradient descent method for minimizing $\mathcal E(X)$ over $(\mathcal S_+^{n,p}, g_{BW})$. 
  With the simple retraction \eqref{eqn:retraction_quotient}, the Riemannian gradient descent method for minimizing the function $h[\pi(Y)]$ on $\mathbb{R}^{n\times p}_*/\mathcal{O}_p$ is given by
  \[Y_{k+1}=Y_{k}-\Delta t 2\nabla \mathcal E(Y_kY_k^T)Y_k,\]
which is the simple Burer-Monteiro gradient descent method
for minimizing $\mathcal E(X)$ over $\mathcal S_+^{n,p}$. See Section 5.1 in \cite{zheng2022riemannian} for details.

\subsubsection{A simple Riemannian Langevin  Monte Carlo scheme}
With the Euler-Maruyama discretization for SDE \eqref{SDE-BW}, and the simple retraction and Riemannian gradient as given previously, 
a simple Riemannian Langevin Monte Carlo scheme for approximating the Riemannian SDE \eqref{SDE-BW} on the Riemannian manifold $({S}^{n,p}_+,g_{BW})$ can be given as
\begin{equation}
    Y_{k+1} = Y_{k} - \Delta t2 \nabla \mathcal{E} (Y_kY_k^T) Y_k  +  \sqrt{\frac{2\Delta t}{\beta}} B_k+\frac{\Delta t}{\beta} U \bmat{\sum_{j: j\neq i}\frac{\sigma_i}{\sigma_i^2 + \sigma_j^2}}_{ii} V^T, 
    \label{Euler_Maruyama_scheme_quotient2} 
\end{equation}
where $B_k$ is $n$-by-$p$ matrix with i.i.d.{} $\mathcal N(0,1)$ entries and
$Y_k = U \Sigma V^T$ is the compact SVD of $Y$ with singular values $\sigma_i>0$ for $i=1,2,\cdots, p$.

Notice that all operations are performed in the space of size $n\times p.$ For finding compact SVD of $Y$, one can first compute QR decomposition of $Y$, which costs $\mathcal O(np^2)+\mathcal O(p^3)$. Then compute SVD of size $p\times p$, which is $\mathcal O(p^3)$. So the complexity of this scheme is $\mathcal O(np^2)+\mathcal O(p^3)$ for each iteration. For large $n$ and small $p$, Scheme BW should be cheaper than Scheme E in each iteration, but they generate different samples for different Gibbs distributions which depend on the metric, i.e., Scheme BW cannot replace Scheme E for generating Gibbs distribution defined by embedded geometry.

\section{Examples with analytical formulae}
\label{sec:example}

In this section, we provide a few examples with analytical formulae so that they can be used in numerical experiments for testing the two schemes \eqref{Euler_Maruyama_scheme_embedded2} and \eqref{Euler_Maruyama_scheme_quotient2} on the Gibbs distribution. 

For the rest of this section, $X=Q\Lambda Q^T\in \mathcal{S}^{n,p}_+$ denotes the full SVD with descending eigenvalues $\lambda_1\geq \lambda_2\geq \cdots \geq \lambda_p>0.$

\subsection{Scalar random variables}

Let $X$ be a random variable satisfying the Gibbs distribution on $\mathcal S^{n,p}_+$ with dimension $N=np-\frac{p(p-1)}2$  under either metric, 
then $X$ is a matrix-valued random variable.
For convenience, we consider a scalar random variable $D=D(X)$ as a function of $X\in \mathcal S^{n,p}_+$, e.g., $D=\|X\|_F$
where $\|\cdot\|_F$ is the matrix Frobenius norm. 

We consider the distribution function for the scalar random variable $D$:
\begin{equation}
    \mathrm{Pr}[D<d]=\frac1{Z_\beta}\int\limits_{U_d}\ee^{-\beta \mathcal{E}}\ud V, \quad
    Z_\beta= \int\limits_{\mathcal M}\ee^{-\beta \mathcal{E}}\ud V,
    \label{scalar-rv}
\end{equation}
 where $U_d:=\{X\in \mathcal{S}^{n,p}_+|D(X)<d\}$ is the domain of integral. For simplicity we only consider symmetric functions such that the random variable $D$, the energy function $\mathcal E$, and the volume form are all invariant under the group action by the orthogonal group $\mathcal O_n$. We consider an energy function $\mathcal E$ satisfying $\mathcal E(X)=\mathcal E(OXO^T)$, $\forall\, O\in\mathcal O_n$, so that Gibbs distribution function only depends on the spectrum of $X$ when considering \eqref{scalar-rv} with $D=||X||_F=\sqrt{\lambda_1^2+\cdots+\lambda_p^2}$. Since $\mathcal O_n$ is an isometry group for both metrics $g_E$ and $g_{BW}$, the volume form $\ud V$ in the two cases is also invariant under $\mathcal O_n$ action.

Notice that $Q$ and $\Lambda$ can be used as coordinates of the manifold $\mathcal S_+^{n,p}$.
The volume form expressed by coordinates $Q$ and $\Lambda$ is given by
\begin{align}
    \ud V=\sqrt{\det g}(\prod_{i=1}^p\ud\lambda_i)\ud \mu_{\mathcal O_n},
\end{align}
where $\mu_{\mathcal O_n}$ is the Haar measure on $\mathcal O_n$, and $g$ is the matrix of metric $g_E$ or $g_{BW}$ expressed under coordinate $Q$ and $\lambda$. For $g_E$ its determinant $\det g$ is 
\begin{align}
    \det g=\big(\prod_{1\leq i<j\leq p}|\lambda_i-\lambda_j|^2\big)\big(\prod_{1\leq i\leq p}\lambda_i^{2(n-p)}\big),
\end{align}
and for $g_{BW}$ it is
\begin{align}
    \det g=\big(\prod_{1\leq i<j\leq p}\frac{|\lambda_i-\lambda_j|^2}{\lambda_i+\lambda_j}\big)\big(\prod_{1\leq i\leq p}\lambda_i^{(n-p)}\big).
\end{align}

So for $g_E$ the distribution $\mathrm{Pr}[D<d]$ is expressed as 
\begin{align}
    \mathrm{Pr}[D<d]=&\frac1{Z_\beta}\int\limits_{||X||_F<d}\ee^{-\beta\mathcal E}\ud V\nonumber\\
    \propto&\int\limits_{\sum\limits_{i=1}^p\lambda_i^2<d^2\atop\lambda_i>0,i=1,...,p}\ee^{-\beta\mathcal E(\lambda_1,...,\lambda_p)}\big(\prod_{1\leq i<j\leq p}|\lambda_i-\lambda_j|\big)\big(\prod_{1\leq i\leq p}\lambda_i^{n-p}\big)\ud \lambda_1\cdots\ud\lambda_p,
\end{align}
where we have used the fact that the integrand does not depend on the coordinate $Q\in\mathcal O_n$, so the integral of $\mu_{\mathcal O_n}$ only provides a constant coefficient. As we could always renormalize $\mathrm{Pr}[D<d]$ by considering the quotient $\frac{\mathrm{Pr}[D<d]}{\mathrm{Pr}[D<\infty]}$, we only need the dependence of the integral on parameter $d$.

Similarly, for the Bures-Wasserstein metric $g_{BW}$ we have 
\begin{align}
    \mathrm{Pr}[D<d]\propto&\int\limits_{\sum\limits_{i=1}^p\lambda_i^2<d^2\atop \lambda_i>0,i=1,...,p}\ee^{-\beta\mathcal E(\lambda_1,...,\lambda_p)}\big(\prod_{1\leq i<j\leq p}\frac{|\lambda_i-\lambda_j|}{\sqrt{\lambda_i+\lambda_j}}\big)\big(\prod_{1\leq i\leq p}\lambda_i^{\frac{n-p}2}\big)\ud \lambda_1\cdots\ud\lambda_p
\end{align}

Next we give a few energy functions.
\subsection{Example I: $\mathcal E(X)=\frac12||X||_F^2$}\label{sec:example1}
This is the simplest example. Using the general expression above, for embedded geometry $g_E$ we have
\begin{align}\label{int1}
    \mathrm{Pr}[D<d]\propto&\int\limits_{\sum\limits_{i=1}^p\lambda_i^2<d^2\atop\lambda_i>0,i=1,...,p}\ee^{-\frac{\beta}2\sum\limits_{i=1}^p\lambda_i^2}\big(\prod_{1\leq i<j\leq p}|\lambda_i-\lambda_j|\big)\big(\prod_{1\leq i\leq p}\lambda_i^{n-p}\big)\ud \lambda_1\cdots\ud\lambda_p\nonumber\\
    =&\int\limits_{0}^d \ee^{-\frac\beta 2 \rho^2}\rho^{N-1} \big(\int\limits_{S^{p-1}_+} 
    \prod_{1\leq i<j\leq p}|\omega_i-\omega_j|\prod_{i=1}^p |\omega_i|^{n-p}\prod_{i=1}^p\ud \omega\big)\ud \rho\nonumber\\
    =&\big(\int\limits_{S^{p-1}_+} 
    \prod_{1\leq i<j\leq p}|\omega_i-\omega_j|\prod_{i=1}^p |\omega_i|^{n-p}\prod_{i=1}^p\ud \omega\big)\int\limits_{0}^d\ee^{-\frac\beta2\rho^2}\rho^{N-1}\ud \rho\nonumber\\
    \propto&\int\limits_0^d \ee^{-\frac\beta2\rho^2}\rho^{N-1}\ud \rho,
\end{align}
where we have used the spherical coordinate for $(\lambda_1,...,\lambda_p)=\rho\omega$, with $\rho=\sqrt{\sum\limits_{i=1}^p\lambda_i^2}$ being the radius and $\omega\in S^{p-1}_+=S^{p-1}\cap\mathbb R^p_+$ being the coordinate on the positive orthant of unit sphere. 

For $g_{BW}$, similarly we have
\begin{align}
    \mathrm{Pr}[D<d]\propto&\int\limits_0^d\ee^{-\beta\rho^2}\rho^{\frac N2-1}\ud \rho.
\end{align}

Now we can see that $\beta D^2=\beta||X||_F^2$ is subject to $\chi^2(N)$ distribution for the embedded metric $g_E$, and $\chi^2(\frac N2)$ distribution for the Bures-Wasserstein metric.

\subsection{Example II: $\mathcal{E}(X)=\mathrm{Tr}(X\log X)$}

We consider the von Neumann entropy $$\mathcal{E}(X)=\mathrm{Tr}(X\log X)=\sum\limits_{i=1}^p\lambda_i\log\lambda_i$$ and construct a more interesting example. The minimizers of $\mathcal E(X)=\mathrm{Tr}(X\log X)$ on $\mathcal{S}^{n,p}_+$ are matrices $X\in \mathcal{S}^{n,p}_+$ with spectrum $\lambda_1=\cdots=\lambda_p=\ee^{-1}$.

The random variable we consider  is still $D=||X||_F$. Since $\mathcal{E}(X)=\mathrm{Tr}(X\log X)=\sum\limits_{i=1}^p\lambda_i\log\lambda_i$ only depends on spectrum, the argument in the previous section about integral on $\mathcal O_n$ still applies. Similar to \eqref{int1}, for $g_E$ we have 
\begin{align*}
    \mathrm{Pr}(D<d)=&\int\limits_{\sum\limits_{i=1}^p\lambda_i^2<d^2\atop \lambda_i>0,i=1,...,p} \ee^{-\beta\sum\limits_{i=1}^p\lambda_i\log\lambda_i}\prod_{1\leq i<j\leq p}|\lambda_i-\lambda_j|\prod_{i=1}^p |\lambda_i|^{n-p}\prod_{i=1}^p\ud \lambda_i    \nonumber\\=&\int\limits_{\sum\limits_{i=1}^p\lambda_i^2<d^2\atop \lambda_i>0,i=1,...,p} \prod_{1\leq i<j\leq p}|\lambda_i-\lambda_j|\prod_{i=1}^p |\lambda_i|^{n-p-\beta\lambda_i}\prod_{i=1}^p\ud \lambda_i,    
\end{align*}
and for $g_{BW}$ we have 
\begin{align*}
    \mathrm{Pr}(D<d)=&\int\limits_{\sum\limits_{i=1}^p\lambda_i^2<d^2\atop \lambda_i>0,i=1,...,p} \prod_{1\leq i<j\leq p}\frac{|\lambda_i-\lambda_j|}{\sqrt{\lambda_i+\lambda_j}}\prod_{i=1}^p |\lambda_i|^{\frac{n-p}2-\beta\lambda_i}\prod_{i=1}^p\ud \lambda_i.    
\end{align*}

Although we do not have a closed expression for both cases, such integrals can be easily approximated by  an accurate quadrature.

\subsection{Example III: $\mathcal{E}(X) =\frac12 ||X-A||_F^2$}
We consider a quadratic function $\mathcal{E}(X) =\frac12 ||X-A||_F^2$ where
 $A\in \mathcal{S}^{n,p}_+$ with $D=||X-A||_F$.
In this example, $\mathcal O_n$ symmetry does not hold, and we can only make an estimate of the distribution function.

The random variable $D$ we are considering now is $D=||X-A||_F$, its distribution function is evaluated as
\begin{align}
    \mathrm{Pr}(D<d)\propto\int_{U_d} \ee^{-\frac\beta2D^2}\ud V
\end{align}
where $U_d=\{X\in\mathcal S^{n,p}_+|D(X)<d\}$. Using delta function, formally we can simplify the integral to 
\begin{align}
    \mathrm{Pr}(D<d)\propto&\int_{\mathcal M}\bd 1_{\{D<d\}}\ee^{-\frac\beta 2D^2}\ud V\\
    =&\int_{\mathcal M}(\int_0^\infty \bd 1_{\{\rho<d\}}\ee^{-\frac\beta2 \rho^2}\delta(D-\rho)\ud \rho)\ud V\nonumber\\
    =&\int_0^\infty\bd 1_{\{\rho<d\}}\ee^{-\frac\beta2 \rho^2}(\int_{\mathcal M} \delta(D-\rho)\ud V)\ud \rho\nonumber\\
    =&\int_0^d \ee^{-\frac\beta2\rho^2}(\int_{\mathcal M}\frac{\ud}{\ud \rho}\bd 1_{\{D-\rho\}}\ud V)\ud\rho\nonumber\\
    =&\int_0^d \ee^{-\frac\beta2\rho^2}\frac{\ud}{\ud \rho}(\int_{\mathcal M}\bd 1_{\{D-\rho\}}\ud V)\ud\rho\nonumber\\
    =&\int_0^d \ee^{-\frac\beta2\rho^2}\frac{\ud}{\ud \rho}V_D(\rho)\ud\rho\nonumber\\
\end{align}
where $V_D(\rho)=\int_{\mathcal M}\bd 1_{\{D<\rho\}}\ud V=\int_{D<\rho}\ud V$.
 
In general it is difficult to calculate $\int_{D<\rho} \ud V$, but we consider the following approximation.
Consider the volume of the ball $B^{n,p}_A(r)=B_A(r)\cap \mathcal{S}^{n,p}_+$, where
$$ B_A(r)=\big\{X\in\mathcal{S}^{n\times n} : ||X-A||_F<r\big\}.$$
     It is difficult to compute $\mathrm{Vol}(B^{n,p}_A(r))$, but we propose the following estimate, for fixed $A\in\mathcal S^{n,p}_+$:
\begin{align}
\label{example3-form}
    \mathrm{Vol}(B^{n,p}_{cA}(r))\approx \alpha r^N,\quad c\gg 1,
\end{align}
where $\alpha$ is a constant that does not depend on $r$, $N$ is the dimension of $\mathcal{S}^{n,p}_+$. For $g_E$, $\alpha$ is exactly the volume of unit ball in $\mathbbm R^N$, while for $g_{BW}$, $\alpha$ depends on dimension $N$ and $A\in\mathcal S^{n,p}_+$.

For the embedded geometry, the  approximation \eqref{example3-form} can be justified by  the following arguments:
\begin{enumerate}
    \item The second fundamental form $\bd{I\!I}_{cA}$ of the manifold is vanishing for fixed $A$ and $c\to\infty$. See \cite{YZLMZ-theory}.
    \item The Riemannian curvature tensor of ambient space $\mathcal S^{n\times n}$ is $0$. Applying the Gauss equation \cite[Prop 3.1]{do1992riemannian} we can express the Riemannian curvature tensor $R$ of $(\mathcal S^{n,p}_+,g_E)$ in terms of its second fundamental form $\bd{I\!I}$:
    \begin{align}
        \langle R(\bd x,\bd y)\bd z,\bd w\rangle=-\langle\bd{I\!I}(\bd x,\bd z),\bd{I\!I}(\bd y,\bd w)\rangle+\langle\bd{I\!I}(\bd x,\bd w),\bd{I\!I}(\bd y,\bd z)\rangle,\\
        \bd x,\bd y,\bd z,\bd w\in T\mathcal S^{n,p}_+,\qquad \langle U,V\rangle=\mathrm{Tr}(UV^T)\text{ is the metric in }\mathcal S^{n\times n},\nonumber.
    \end{align}
    Thus, with vanishing $\bd{I\!I}$ we have vanishing Riemannian curvature tensor, and zero sectional curvature.
    \item Vanishing extrinsic curvature and intrinsic curvature means that the neighborhood is approximately an Euclidean space, so the ball $B^{n,p}_{cA}(r)$ is approximately just a ball in $\mathbbm R^N$ and has volume $\alpha r^N$, with $\alpha$ being the volume of a unit ball.
\end{enumerate}

For the $g_{BW}$ metric, following similar arugments, we can get the same approximation \eqref{example3-form}. We emphasize that the approximation \eqref{example3-form} is accurate only if $c$ is large enough. Putting all this together, when $A$ has eigenvalues $\lambda_1\geq \cdots\geq \lambda_p\gg 1$, we have the following 
\begin{equation}
 \mathrm{Pr}(D<d)\propto\int\limits_{D<d} \ee^{-\frac\beta2D^2}\ud V=\int\limits_0^d \ee^{-\frac\beta2\rho^2} \frac{\ud}{\ud \rho}\big( \int\limits_{D<\rho} \ud V\big)\ud \rho \appropto \int\limits_0^t \ee^{-\frac\beta2\rho^2} \rho ^{N-1} \ud \rho,
\end{equation}
 where $\appropto$ stands for {\it being approximately proportional to}.

\subsection{MCMC numerical integration}
\label{sec: quad}
It is well known that MCMC  can be used for integrating a function numerically, and that one of the main advantages is that the convergence rate is independent of the dimension. 
Both schemes in this paper are MCMC type sampling schemes on the manifold. 
Suppose we have generated samples $X_i$ satisfying the Gibbs distribution on the manifold, e.g., 
$$X_i\sim \frac{1}{Z_\beta} e^{-\beta \mathcal E(X)} \ud V_g,$$
where $Z_\beta=\int\limits_{\mathcal \mathcal{S}^{n,p}_+}\ee^{-\beta \mathcal E(X)}\ud V$ is an unknown normalization factor and $\ud V$ is the volume form depending on the metric. 
Then for approximating the integral of a nice function $f(X)$ on the same manifold
$ \int\limits_{\mathcal \mathcal{S}^{n,p}_+}f(X)\ud V,$ we can use  
\begin{align}
    \frac1m\sum_{i=1}^m f(X_i)\ee^{\beta \mathcal E(X_i)}\approx \frac{\int\limits_{\mathcal \mathcal{S}^{n,p}_+} f(X)\ud V}{\int\limits_{\mathcal \mathcal{S}^{n,p}_+}\ee^{-\beta \mathcal E(X)}\ud V}=\frac{1}{Z_\beta} \int\limits_{\mathcal \mathcal{S}^{n,p}_+} f(X)\ud V,
    \label{mc-quadrature}
\end{align}
because each $f(X_i)\ee^{\beta \mathcal E(X_i)}$ is a random variable with expectation 
$$\mathbb E\left[f(X_i)\ee^{\beta \mathcal E(X_i)}\right]=\frac{1}{Z_\beta} \int\limits_{\mathcal \mathcal{S}^{n,p}_+} f(X_i)\ee^{\beta \mathcal E(X_i)} \ee^{-\beta \mathcal E(X_i)}\ud V,$$
and the left hand side is a random variable with expectation 
\[\mathbb E\left[\frac1m\sum_{i=1}^m f(X_i)\ee^{\beta \mathcal \mathbb E(X_i)}\right]=\frac1m\sum_{i=1}^m \mathbb E\left[ f(X_i)\ee^{\beta \mathcal \mathbb E(X_i)}\right]=\frac{1}{Z_\beta} \int\limits_{\mathcal \mathcal{S}^{n,p}_+} f(X)\ud V,  \]
where the expectation $\mathbb E[\cdot]$ is taken w.r.t. Gibbs distribution under corresponding metric.

So using the generated samples $X_i$, we can approximate the integral $\int\limits_{\mathcal{S}^{n,p}_+}f(X)\ud V$ up to a constant  $Z_\beta$ that does not depend on $f(X)$. Notice that the additional advantage of Monte Carlo type quadrature on a  manifold  is that we do not need to know what $\ud V$ is. On the other hand, $Z_\beta$ cannot be approximated by the same approach.  
Though we do not consider any specific application for numerical integration, equation \eqref{mc-quadrature} can be used as one way to validate the Riemannian Langevin Monte Carlo schemes. 

For the following special functions, it is possible to calculate exact integrals.
For the energy function $\mathcal{E}(X) = \frac{1}{2}\norm{X}^2_F$,
and a special integrand $f(X) = \norm{X}^k_F \ee^{- \frac{\alpha}{m} \norm{X}_F^m}$ with $k>-N, m>2, \alpha>0$, using the results in \ref{sec:example1}, the distribution of $D=||X||_F$ is
\begin{align}
    \text{ for metric }g_E:&\quad\mathrm{Pr}[D<d]\propto\int_0^d \ee^{-\frac\beta2\rho^2}\rho^{N-1}\ud\rho,\\
    \text{ for metric }g_{BW}:&\quad\mathrm{Pr}[D<d]\propto\int_0^d \ee^{-\frac\beta2\rho^2}\rho^{\frac N2-1}\ud\rho,
\end{align}
so the integral on the manifold could be expressed by expectation of a random variable, which leads to
 
\begin{align}\label{quad-example-1}
    \text{ for }g_E:&\frac1{Z_\beta}\int_{\mathcal S^{n,p}_+}f(X)\ud V=\mathbb E[f(X)\ee^{\frac\beta2 ||X||_F^2}]=\mathbb E[D^k\ee^{-\frac\alpha m D^m}\ee^{\frac\beta2 D^2}]\nonumber\\
    =&\frac{\int_0^\infty \rho^k\ee^{-\frac\alpha m \rho^m+\frac\beta2\rho^2}\rho^{N-1}\ee^{-\frac\beta2\rho^2}\ud\rho}{\int_0^\infty\rho^{N-1}\ee^{-\frac\beta2\rho^2}\ud\rho}=\frac{\frac{1}{m}(\alpha/m)^{-\frac{k+N}{m}} \Gamma((k+N)/m)}{\frac{1}{2}(\beta/2)^{-N/2}\Gamma(N/2)}\\
    \text{ for }g_{BW}:&\frac1{Z_\beta}\int_{\mathcal S^{n,p}_+}f(X)\ud V=\mathbb E[f(X)\ee^{\frac\beta2 ||X||_F^2}]=\mathbb E[D^k\ee^{-\frac\alpha m D^m}\ee^{\frac\beta2 D^2}]\nonumber\\
    =&\frac{\int_0^\infty \rho^k\ee^{-\frac\alpha m \rho^m+\frac\beta2\rho^2}\rho^{\frac N2-1}\ee^{-\frac\beta2\rho^2}\ud\rho}{\int_0^\infty\rho^{\frac N2-1}\ee^{-\frac\beta2\rho^2}\ud\rho}=\frac{\frac{1}{m}(\alpha/m)^{-\frac{k+N/2}{m}} \Gamma((k+N/2)/m)}{\frac{1}{2}(\beta/2)^{-N/4}\Gamma(N/4)}. \label{quad-example-2}
\end{align}

\section{Numerical tests} 
\label{sec:tests}
In this section we test the samples generated by the two Riemannian Langevin Monte Carlo schemes   \eqref{Euler_Maruyama_scheme_embedded2} and \eqref{Euler_Maruyama_scheme_quotient2} on the examples constructed in the previous section. 
The samples are generated by the following procedure: we run the iterative schemes \eqref{Euler_Maruyama_scheme_embedded2} or \eqref{Euler_Maruyama_scheme_quotient2} for sufficiently many $\tilde m$  iterations then take the last $m$ iterates as the samples for the Gibbs distribution. Both $\tilde m$ and $m$ should be chosen such that the $(\tilde m-m)$-th iterate has already reached equilibrium  e.g., $\tilde m$ is $6,000,000$ and $m$ is $5,000,000$ for specially chosen energy functions and parameters $\beta$. 

Now
suppose we have generated   samples $X_i\in \mathcal S_+^{n,p}$ $(i=1,\cdots, m)$ for either metric. In order to test or show the numerical convergence to the Gibbs distribution, we will consider two kinds of numerical tests.

The first kind of tests is to test on the scalar random variable $D(X)=\|X\|_F$ or $D(X)=\|X-A\|_F$ as described in Section \ref{sec:example}.  
 Then we compare the cumulative distribution function (CDF) of the random variable $D$ with its empirical CDF calculated from the MCMC samples. 

Denote the true CDF of $D$ by $F_D(t):= \Pr (D \leq t)$.
The empirical CDF of samples is
$$\hat{F}_D(t) := \frac{1}{m} \sum_{i=1}^m \mathbbm{1}_{D(X_i) \leq t},
$$  
where $\mathbbm{1}_{D(X_i) \leq t}$ takes value $1$ if $D(X_i) \leq t$,  and value $0$ if otherwise. 
The Kolmogorov–Smirnov test statistic (K-S statistic)
  is defined by 
\begin{equation}\label{eqn:ks}
    KS_{D} :=  \sup_t\abs{F_D(t) - \hat{F}_D(t)}. 
\end{equation}

In our numerical tests, we compute the KS statistic by taking the maximum difference of $F_D$ and $\hat{F}_D$ at 100 equally spaced points in the interval $[0,t_{max}]$ where $F_D(t_{max}) \approx 1$.

The second kind of tests is 
on the integral examples in Section \ref{sec: quad}, let $X$ be a random variable satisfying Gibbs distribution on the manifold $\mathcal S_+^{n,p}$ under either metric.
Define 
\[
\mu := \E{f(X) \ee ^{\beta \mathcal{E}(X)}}=\frac{1}{Z_\beta} \int\limits_{\mathcal{S}^{n,p}_+} f(X) \ud V.
\] 
Given $m$ samples $X_i\in \mathcal S_+^{n,p}$, we define
\begin{equation}\label{mu-hat}
\hat \mu_m := \frac{1}{m} \sum_{i=1}^m f(X_i) \ee ^{\beta \mathcal{E} (X_i) }.
\end{equation}
Notice that samples generated by MCMC are not independent.
If we assume
\[
\sigma^2 := \mbox{var}\left( f(X_1) \ee ^{\beta \mathcal{E}(X_1)}\right)  +2 \sum_{k=1}^\infty \mbox{cov}\left( f(X_1) \ee ^{\beta \mathcal{E}(X_1)}, f(X_{1+k}) \ee ^{\beta \mathcal{E}(X_{1+k})} \right)   < \infty,
\]

then by the Markov Chain Central Limit Theorem\cite{jones2004markov, geyer1998markov}, as $m\to \infty$,  we have 
\begin{equation}\label{eqn:mcclt}
\sqrt{m}(\hat{\mu}_m - \mu) \rightarrow \mathcal{N}(0,\sigma^2)
\end{equation}
where the convergence is in the sense of distribution. 
Thus if $m\gg 1$, $\frac{\hat{\mu}_m - \mu}{\mu}$  roughly follows the distribution $\mathcal{N}(0,\mathcal O(\frac1m))$ and the relative error term $\abs{\frac{\hat{\mu}_m - \mu}{\mu}}$ roughly follows the folded normal distribution with mean $\mathcal O(\frac{1}{\sqrt{m}} )$ and variance $\mathcal O(\frac1m)$.  Hence we can use $\hat{\mu}_m$ defined in \eqref{mu-hat} to estimate $\mu = \frac{1}{Z_\beta}  \int\limits_{\mathcal{S}^{n,p}_+} f(X) \ud V$,
and the relative error is $\mathcal O(\frac{1}{\sqrt{m}} )$.

\subsection{Numerical validation of the scalar variable $D(X)$}

The manifold $\mathcal S_+^{n,p}$ has dimension $N = np - p(p-1)/2$.
For both metrics, we consider three examples in Section \ref{sec:example} with special energy functions $\mathcal{E}$ in the Gibbs distribution $\ee ^{-\beta \mathcal{E}}$ and the CDF for the scalar variable $D(X)$: 
\begin{enumerate}
    \item  Example I:  $\mathcal{E}(X) = \frac{1}{2} \norm{X}_F^2$  with the CDF for $D(X)=\|X\|_F$:
    \begin{equation*}
  \mbox{For $g_E$}:\quad  F_D(t) = \Pr(\norm{X}_F \leq t) \propto \int\limits_0^t \ee^{-\frac\beta2\rho^2}\rho^{N-1}\ud \rho,
\end{equation*}
\begin{equation*}
   \mbox{For $g_{BW}$}:\quad   F_D(t) = \Pr(\norm{X}_F \leq t) \propto \int\limits_0^t \ee^{-\frac\beta2\rho^2}\rho^{N/2-1}\ud \rho.
\end{equation*}
    \item Example II: $\mathcal{E}(X) = \Tr(X\log X)$ 
     with the CDF $F_D(t) =  \Pr(\norm{X}_F \leq t)$ for $D(X)=\|X\|_F$:
    \begin{equation*}
 \mbox{For $g_E$}:\quad   F_D(t)   \propto \int\limits_{\sum\limits_{i=1}^p\lambda_i^2<t^2\atop \lambda_i>0,i=1,...,p} \prod_{1\leq i<j\leq p}|\lambda_i-\lambda_j|\prod_{i=1}^p |\lambda_i|^{n-p-\beta\lambda_i}\prod_{i=1}^p\ud \lambda_i,
\end{equation*}
\begin{equation*}
\mbox{For $g_{BW}$}:\quad     F_D(t)   \propto \int\limits_{\sum\limits_{i=1}^p\lambda_i^2<t^2\atop \lambda_i>0,i=1,...,p} \prod_{1 \leq i < j \leq p} \frac{|\lambda_i - \lambda_j|}{\sqrt{\lambda_i + \lambda_j}}\prod_{i=1}^p |\lambda_i|^{\frac{n-p-1}{2}-\beta\lambda_i}\prod_{i=1}^p\ud \lambda_i.
\end{equation*}
which is a $p$-fold integral and can be approximated  accurately by quadrature such as Simpson's rule for relatively small values of $p$, e.g., $p=2,3$. 
    \item Example III: $\mathcal{E}(X) = \frac{1}{2} \norm{X-A}_F^2$ where $A\in\mathcal S_+^{n,p}$ has eigenvalues $\lambda_1\geq \cdots\geq \lambda_p\gg 1$,
     with the CDF for $D(X)=\|X-A\|_F$:
    \begin{equation*}
 \mbox{For both $g_E$ and $g_{BW}$}:\quad   F_D(t) = \Pr(\norm{X-A}_F \leq t) \appropto \int\limits_0^t \ee^{-\frac\beta2\rho^2} \rho ^{N-1} \ud \rho.
\end{equation*} 
\end{enumerate}

\begin{figure}[htbp]
\centering
\subfigure[Scheme E \eqref{Euler_Maruyama_scheme_embedded2} on $(\mathcal S^{n,p}_+,g_E)$
with $\Delta t=0.001$ and $\beta=0.4$. The error between two CDFs is  $KS= 0.0054$.]{
\includegraphics[width=0.4\textwidth]{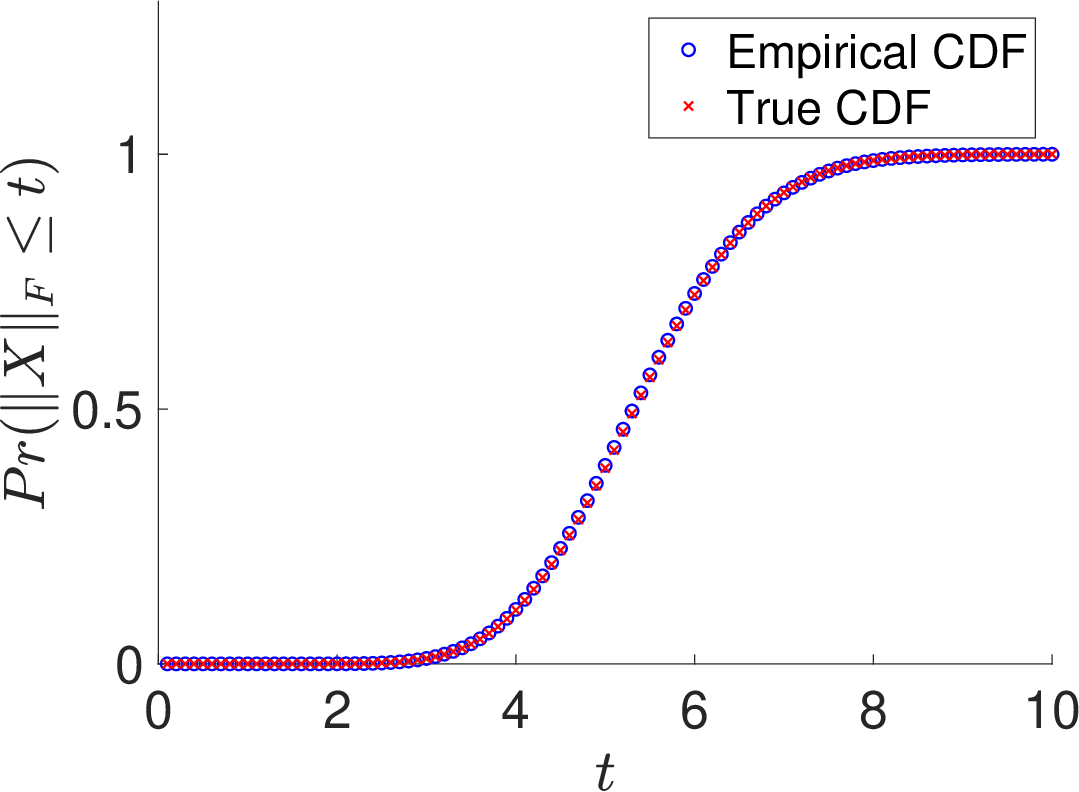}
}
\hspace{0.6cm}
\subfigure[Scheme BW \eqref{Euler_Maruyama_scheme_quotient2} 
 on $(\mathcal S^{n,p}_+,g_{BW})$ with $\Delta t=0.001$ and $\beta=0.4$. The error between two CDFs is  $KS=0.0023$.]{
\includegraphics[width=0.4\textwidth]{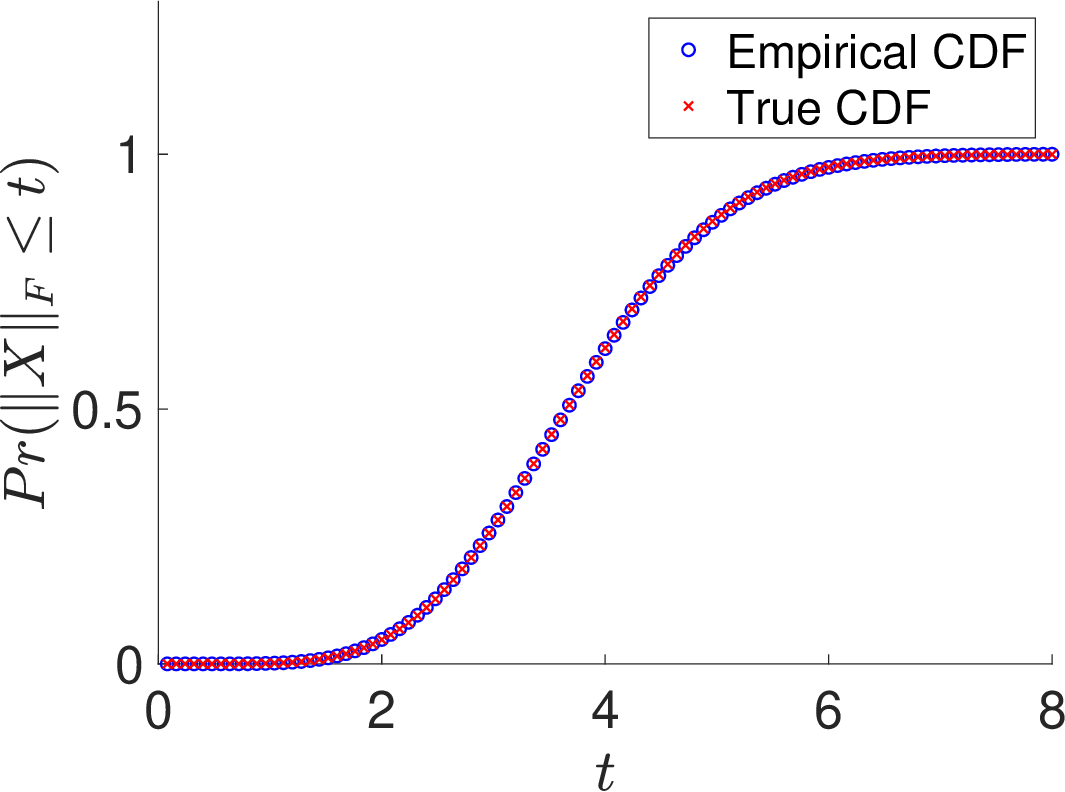}
}
\caption{Example I: $\mathcal{E}(X) = \frac{1}{2}\norm{X}_F^2, n=5,p=3$ and manifold dimension is $N=12$. The empirical CDF is computed by $5E6$ MCMC samples generated after $6E6$ iterations of the Riemannian Langevin Monte Carlo schemes. Both CDFs of scheme E and scheme BW are evaluated at 100 equally spaced points on $[0,10]$ and $[0,8]$, respectively, and the difference can be measured by the KS statistic \eqref{eqn:ks}. }
\label{fig:example-1}
\end{figure}

\begin{figure}[htbp]
    \centering
    \subfigure[
    Scheme E \eqref{Euler_Maruyama_scheme_embedded2} on $(\mathcal S^{n,p}_+,g_E)$
with $\Delta t=0.001$ and $\beta=0.5$. The error between two CDFs is  $KS= 0.0096$
]{
    \includegraphics[width=0.4\textwidth]{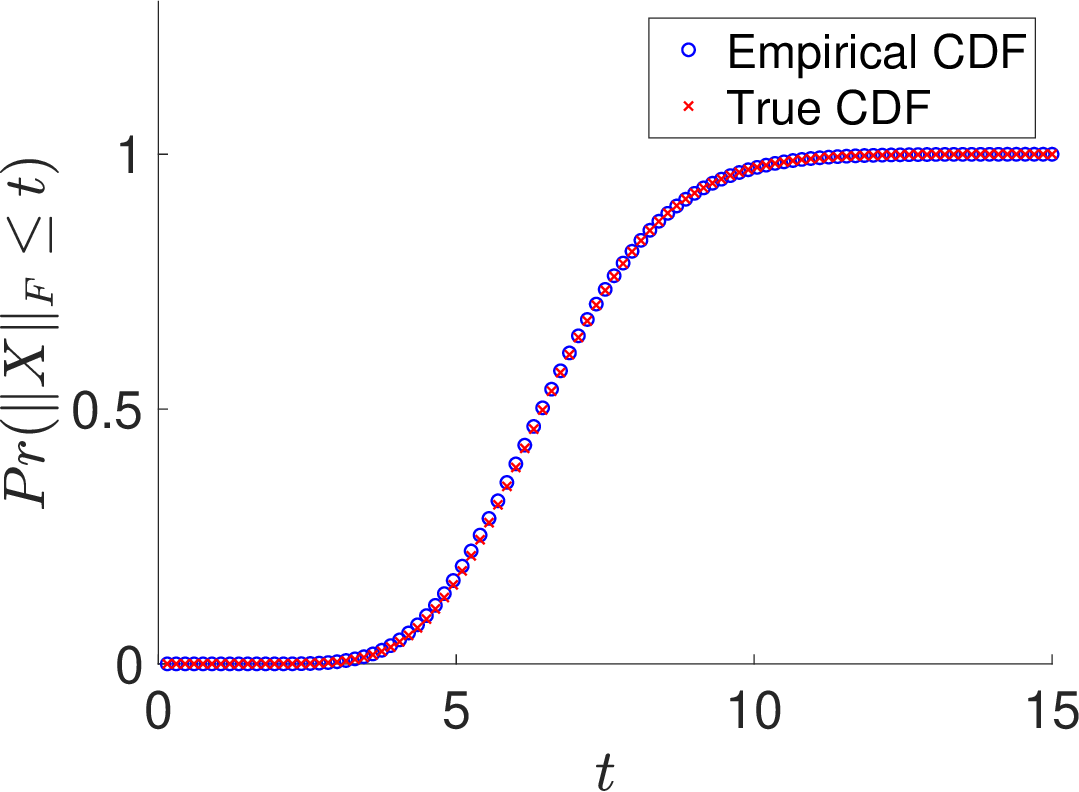}
    }
    \hspace{0.6cm}
    \subfigure[Scheme BW \eqref{Euler_Maruyama_scheme_quotient2} 
 on $(\mathcal S^{n,p}_+,g_{BW})$ with $\Delta t=0.001$ and $\beta=0.5$. The error between two CDFs is $KS= 0.0043$. ]{
    \includegraphics[width=0.4\textwidth]{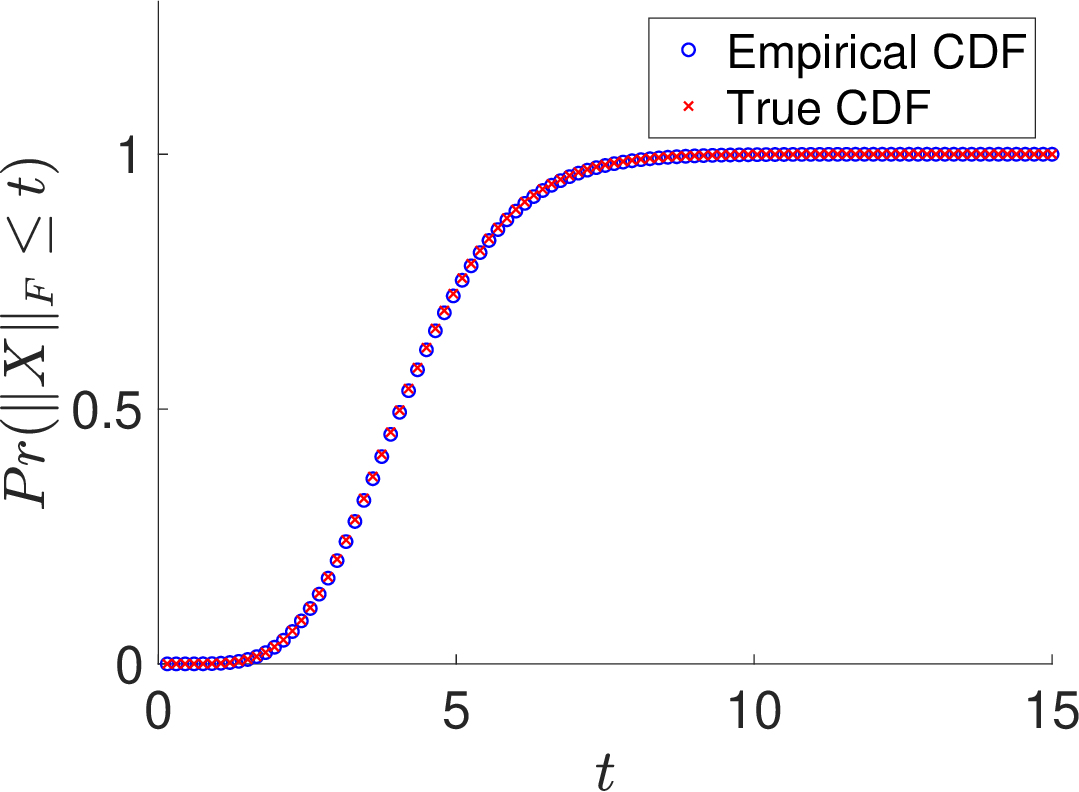}
    }
    \caption{Example II: $\mathcal{E}(X) = \Tr(X \log X)$,  $n=5, p=3$ and manifold dimension is $N=12$. The empirical CDF is computed by $5E6$ MCMC samples generated after $6E6$ iterations of the Riemannian Langevin Monte Carlo schemes. Both CDFs are evaluated at 100 equally spaced points on $[0,15]$, and the difference can be measured by the KS statistic \eqref{eqn:ks}. }
    \label{fig:example-2-1}
\end{figure}

\begin{figure}[htbp]
    \centering
    \subfigure[Scheme E \eqref{Euler_Maruyama_scheme_embedded2} on $(\mathcal S^{n,p}_+,g_E)$
with $\Delta t=0.001$ and $\beta=0.5$. The error between two CDFs is  $KS =0.006$.]{
    \includegraphics[width=0.4\textwidth]{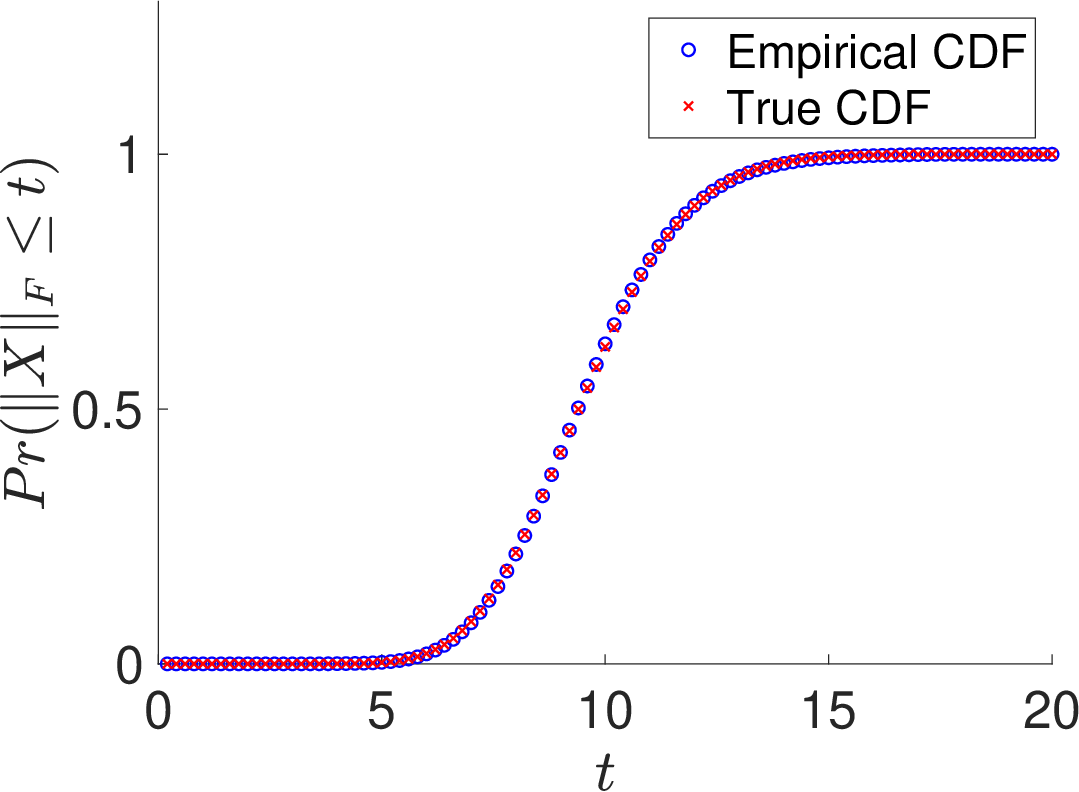}
    }
    \hspace{0.6cm}
    \subfigure[Scheme BW \eqref{Euler_Maruyama_scheme_quotient2} 
 on $(\mathcal S^{n,p}_+,g_{BW})$ with $\Delta t=0.001$ and $\beta=0.5$. The error between two CDFs is $KS=0.0043$.]{
    \includegraphics[width=0.4\textwidth]{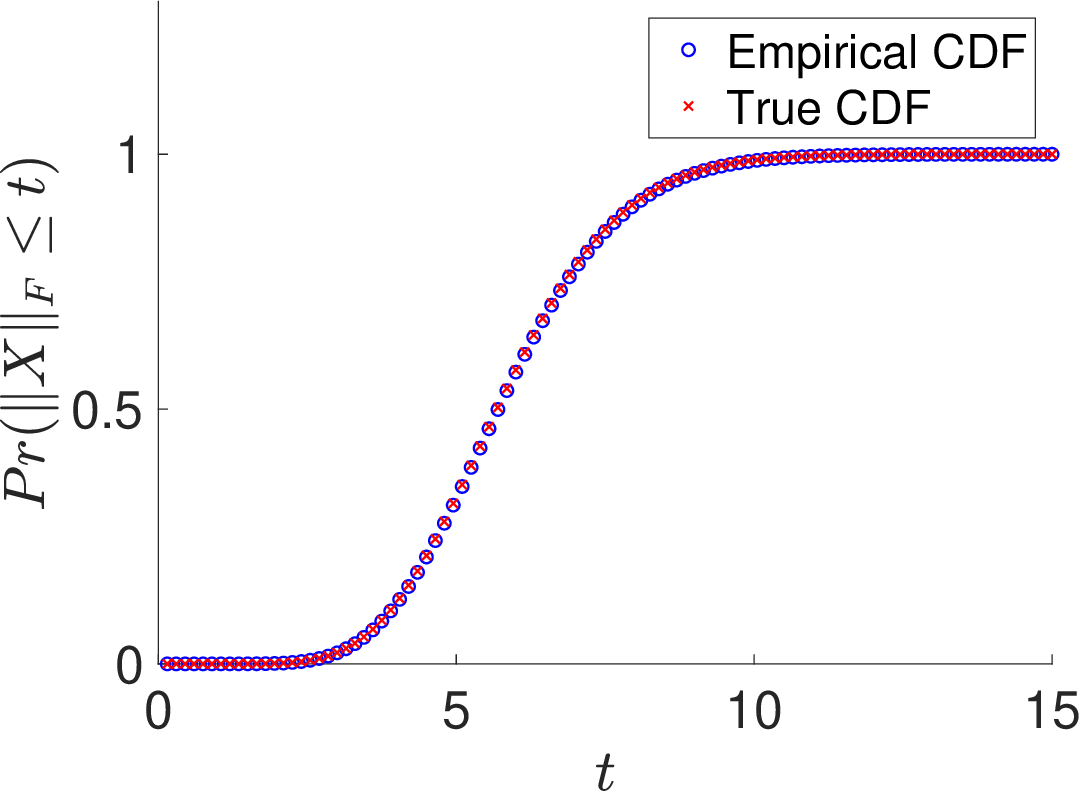}
    }
    \caption{Example II: $\mathcal{E}(X) = \Tr(X \log X)$,  $n=10, p=2$ and manifold dimension is $N=19$. The empirical CDF is computed by $5E6$ MCMC samples generated after $6E6$ iterations of the Riemannian Langevin Monte Carlo schemes. Both CDFs of scheme E and scheme BW are evaluated at 100 equally spaced points on $[0,20]$ and $[0,15]$,respectively, and the difference can be measured by the KS statistic \eqref{eqn:ks}. }
    \label{fig:example-2-2}
\end{figure}

\begin{figure}[htbp]
    \centering
    \subfigure[Scheme E \eqref{Euler_Maruyama_scheme_embedded2} on $(\mathcal S^{n,p}_+,g_E)$
with $\Delta t=0.001$ and $\beta=0.4$. The error between two CDFs is  $KS = 0.0084$.]{
    \includegraphics[width=0.4\textwidth]{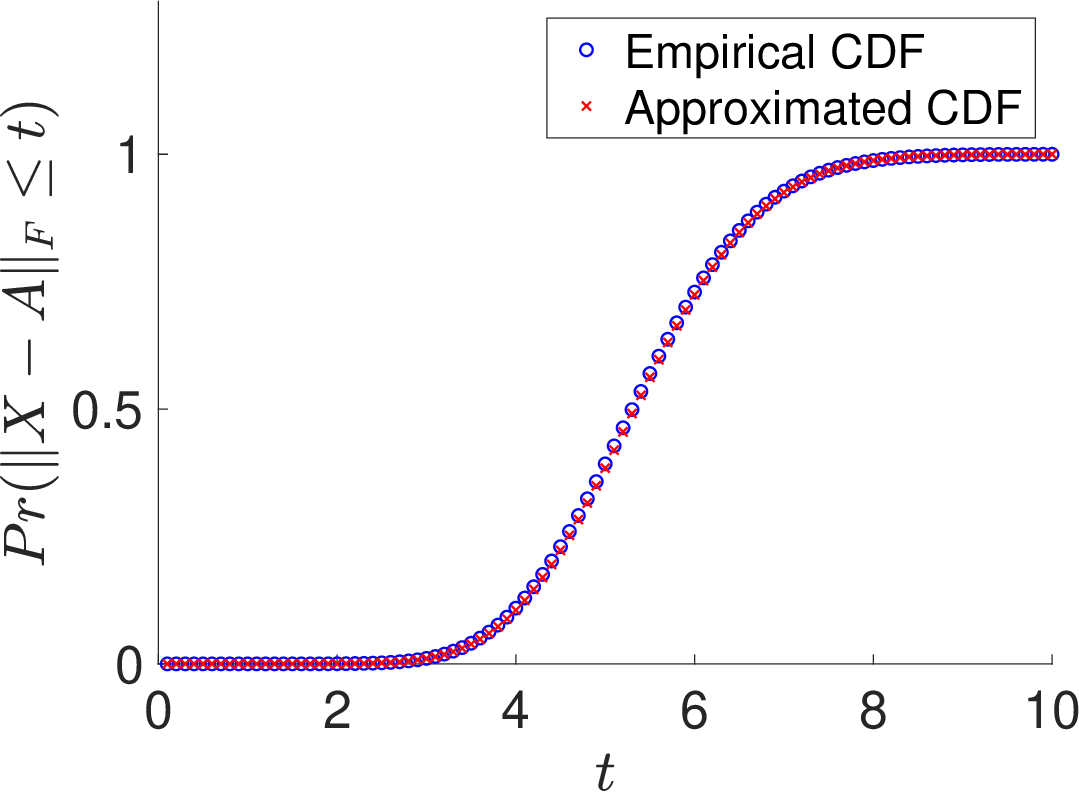}
    }
    \hspace{0.6cm}
    \subfigure[Scheme BW \eqref{Euler_Maruyama_scheme_quotient2} 
 on $(\mathcal S^{n,p}_+,g_{BW})$ with $\Delta t=$2E-7 and $\beta=0.4$. The error between two CDFs is $KS=0.0052$.]{
    \includegraphics[width=0.4\textwidth]{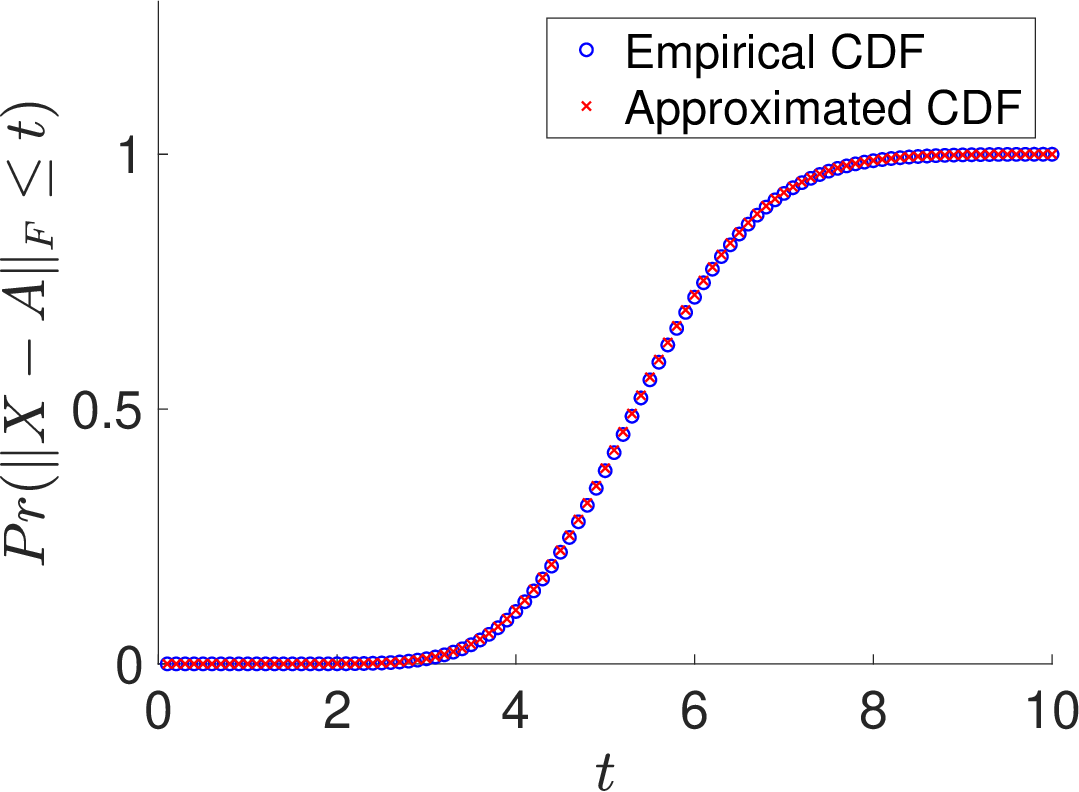}
    }
    \caption{Example III: $\mathcal{E}(X) = \frac{1}{2}\norm{X-A}_F^2$, $n=5, p=3$ and manifold dimension is $N=12$. The nonzero eigenvalues of $A$ are equally spaced between 10000 and 20000.  The empirical CDF is computed by $5E6$ MCMC samples generated after $6E6$ iterations of the Riemannian Langevin Monte Carlo schemes. Both CDFs are evaluated at 100 equally spaced points on $[0,10]$, and the difference can be measured by the KS statistic \eqref{eqn:ks}. }
    \label{fig:example-3}
\end{figure}

In implementation of the scheme, the step size $\Delta t$ and $\beta$ in schemes  \eqref{Euler_Maruyama_scheme_embedded2} and \eqref{Euler_Maruyama_scheme_quotient2} are two parameters that need to be tuned
to reach equilibrium with reasonable computing time. 
We first use a numerically stable $\Delta t$ then adjust $\beta$ so that the noise term has reasonable variance. And of course one needs a sufficient large number of iterations for schemes \eqref{Euler_Maruyama_scheme_embedded2} and \eqref{Euler_Maruyama_scheme_quotient2}  to reach their equilibrium state, and a sufficient large number $m$ of samples to observe numerical convergence toward the Gibbs distribution through the scalar random variable $D$, e.g., the KS statistic \eqref{eqn:ks} should be small.  See Figure \ref{fig:example-1},  Figure \ref{fig:example-2-1}, Figure \ref{fig:example-2-2}, and Figure \ref{fig:example-3} for the numerical results. 
   
\subsection{MCMC numerical integration}
We consider special cases $k=0, m=2$ in the examples \eqref{quad-example-1} and \eqref{quad-example-2}, then  \eqref{quad-example-1} reduces to $(\frac{\beta}{\alpha})^{N/2}$
and \eqref{quad-example-2} reduces to $(\frac{\beta}{\alpha})^{N/4}$.
In other words, we may verify the numerical convergence of samples $X_i$ to Gibbs distribution by verifying
\begin{equation}\label{eqn:convergence_sample_mean_eg1-2}
 \mbox{For $g_{E}$}:\quad  \frac{1}{m} \sum_{i=1}^m \ee ^{-\frac{\alpha - \beta}{2}\norm{X_i}^2_F} \to  (\frac{\beta}{\alpha})^{N/2},
\end{equation}
\begin{equation}\label{eqn:convergence_sample_mean_scheme2-eg1-1}
 \mbox{For $g_{BW}$}:\quad     \frac{1}{m} \sum_{i=1}^m \ee ^{-\frac{\alpha - \beta}{2}\norm{X_i}^2_F} \to  (\frac{\beta}{\alpha})^{N/4}.
\end{equation}

\begin{figure}[htpb]
    \centering
    \subfigure[Integration on $(\mathcal S^{n,p}_+,g_E)$ via samples generated by Scheme E \eqref{Euler_Maruyama_scheme_embedded2} 
with $\Delta t=0.001$ and $\beta=0.4$.]{
    \includegraphics[width=0.45\textwidth]{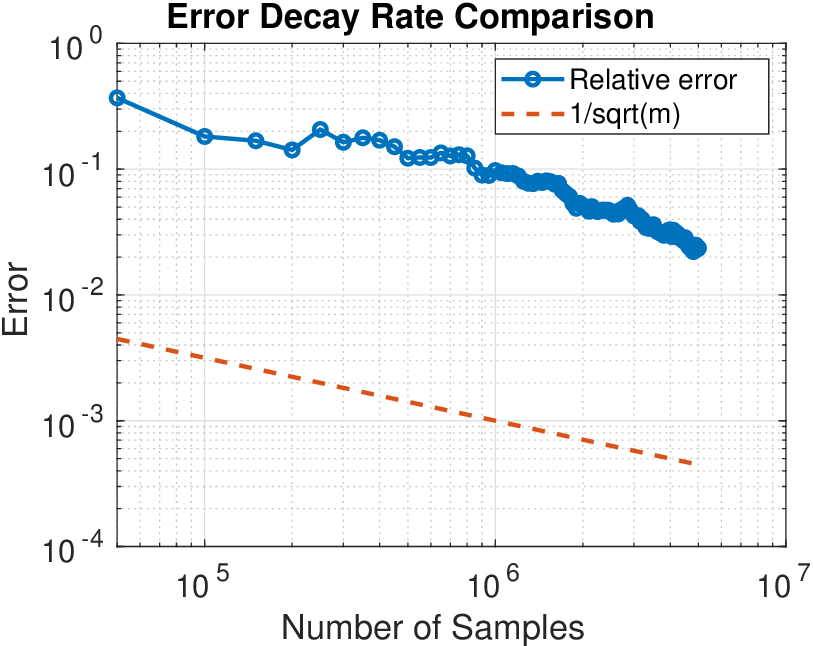}}
    \hspace{0.5cm}
    \subfigure[Integration on $(\mathcal S^{n,p}_+,g_{BW})$ via samples generated by  Scheme BW \eqref{Euler_Maruyama_scheme_quotient2} 
  with $\Delta t=0.001$ and $\beta=0.4$.]{
    \includegraphics[width=0.45\textwidth]{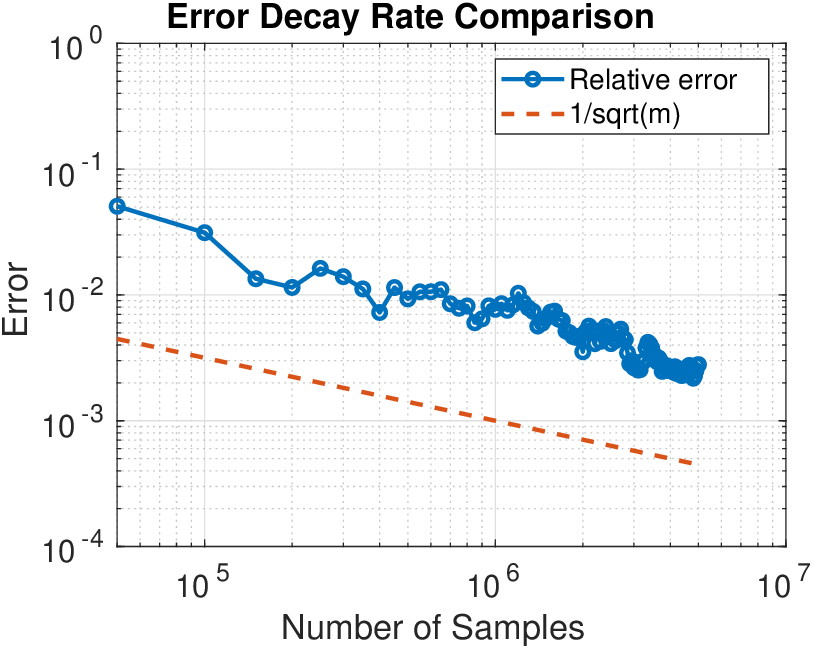}
    }
    \caption{Convergence rate of the relative error of $\abs{\frac{\hat{\mu}_m - \mu}{\mu}}$ MCMC integration on the manifold with $n=10,p=2$ and dimension $N=19$.   
    Parameters are $\alpha = 0.75, \beta = 0.4$, for which it is a numerical integration of the function  $f(X) = \frac{1}{2}\norm{X}_F^2$ on the manifold $ \mathcal S_+^{n,p}$.
    The error shown is the averaged one of 12 independent runs. }
    \label{fig:quad}
\end{figure}

In Figure \ref{fig:quad} we indeed observe the $\mathcal O(1/\sqrt{m})$ for the relative error of numerical integration.

\subsection{A numerical study of the convergence to equilibrium}

\begin{figure}[htbp]
    \centering
\subfigure[Convergence of
    SDEs to the equilibrium.]{\includegraphics[width = 0.8\textwidth]{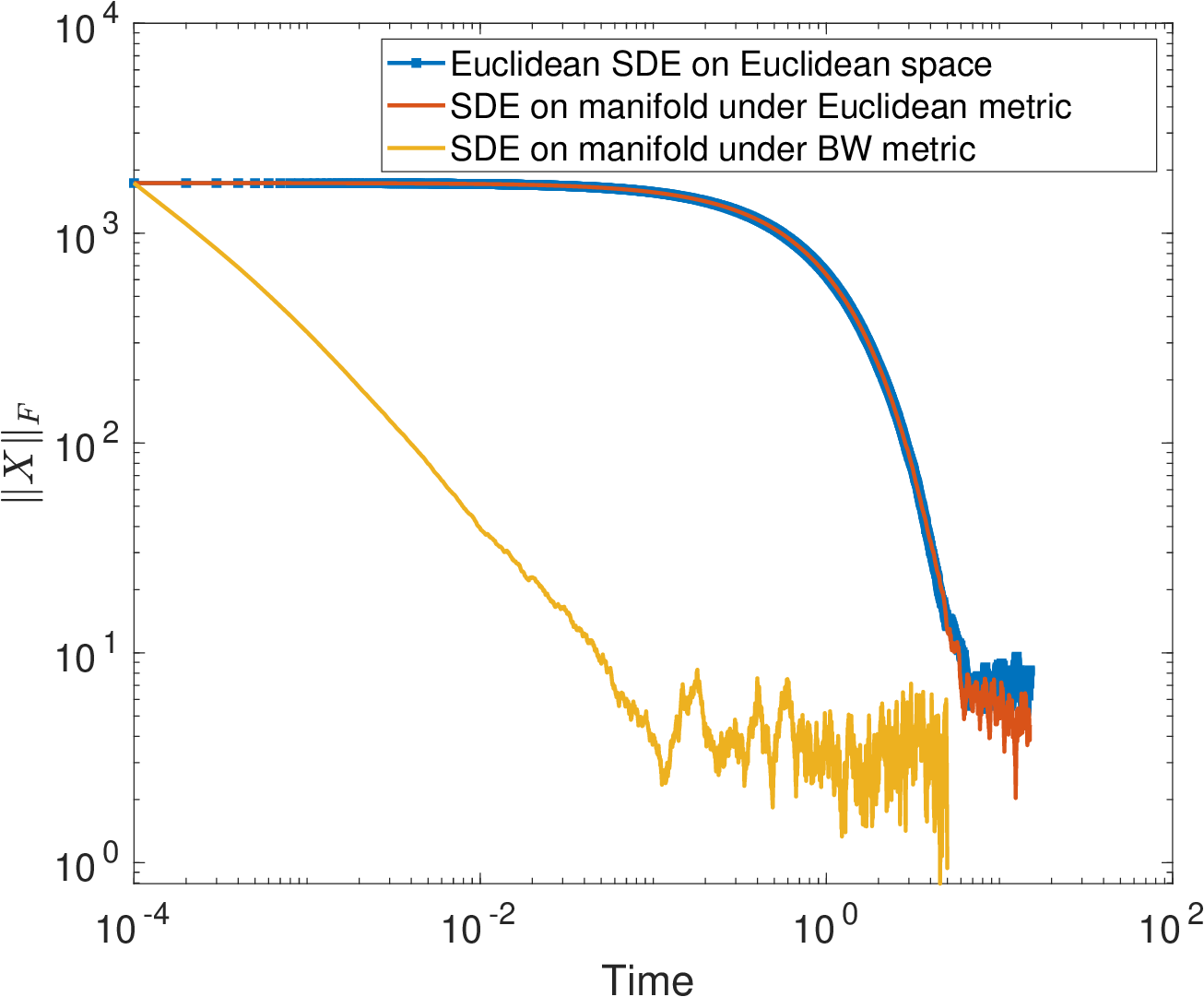}}\\
\subfigure[The true CDF of $D=\norm{X}_F$  on $(\mathcal S_+^{n,p}, g_{E})$. This implies the blue curve in Figure (a) has not reached equilibrium at Time=1.]{\includegraphics[width=0.4\textwidth]{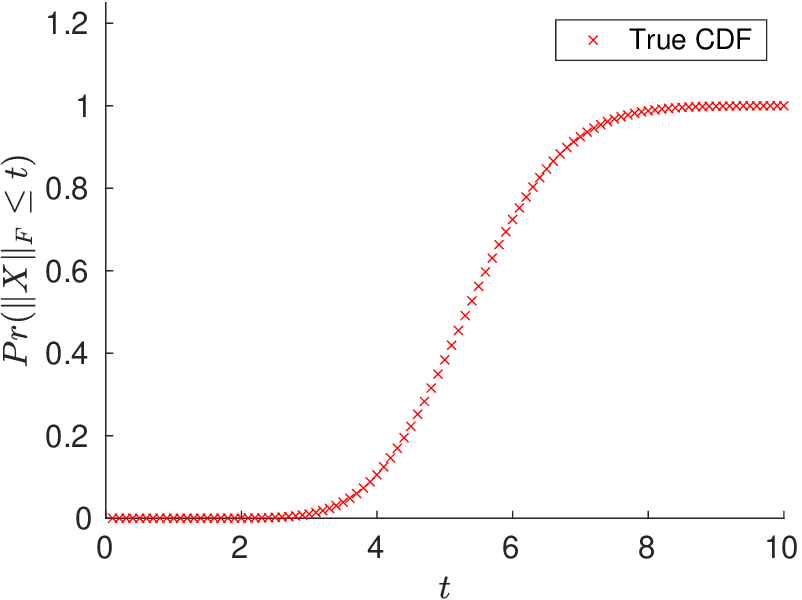}}
\hspace{1cm}
    \subfigure[The true CDF of $D=\norm{X}_F$  on $(\mathcal S_+^{n,p}, g_{BW})$, and empirical CDF 
    of 5E6 samples from the Scheme BW after Time=5.]{  \includegraphics[width=0.4\textwidth]{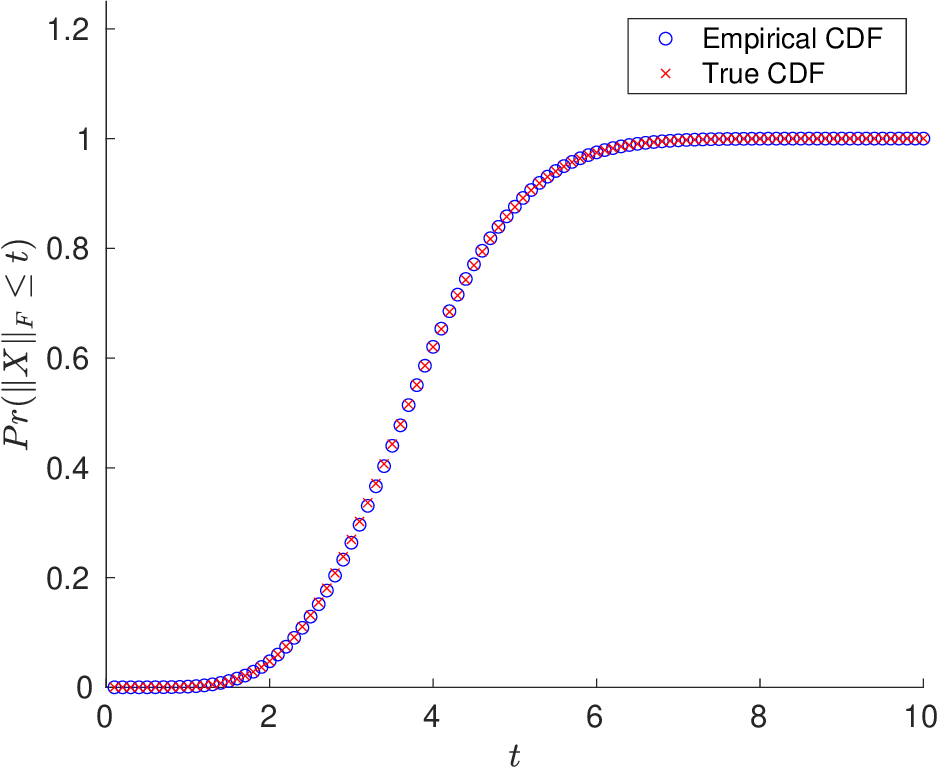}}
    \caption{ For $\mathcal S_+^{n,p}$ with $n=5,p=3$, the dimension is $N=12$, and $\mathbb R^{n\times n}$ has dimension $25.$ The Gibbs measure is $\ee^{-\beta \mathcal E(x)}$ with $\beta=0.4$ and $\mathcal E=\|X\|_F^2$. The time step size is $1e-4$. The initial guess is a random PSD matrix of rank-3 with eigenvalues = [1000,1000,1000].}
    \label{fig:sde-convergence}
\end{figure}

The general mathematical theory of convergence of a Langevin equation to its equilibrium measure has been well studied; we consider the specific case of the RLE studied here in the companion paper~\cite{YZLMZ-theory}. One particular application of the two Riemannian Langevin Monte Carlo schemes is to use them to numerically study the SDE solutions, e.g., by taking very small time steps, a Riemannian Langevin Monte Carlo scheme approximates the Riemannian Langevin equation on the manifold. We have shown comparison of the Langevin equation on $(\mathcal S_+^{n,p}, g_{E})$, $(\mathcal S_+^{n,p}, g_{BW})$, $\mathbb R^{n\times n}$ in Figure \ref{fig:sde-convergence}, in which we can see interesting differences between two metrics.  With all three figures in Figure \ref{fig:sde-convergence}, we can see that the SDE on $(\mathcal S_+^{n,p}, g_{BW})$ has a much faster convergence to its Gibbs measure than the SDE on $(\mathcal S_+^{n,p}, g_{E})$.

\section{Conclusion}\label{sec:remark}

We have constructed two efficient Riemannian Langevin Monte Carlo schemes for sampling PSD matrices of fixed rank from the Gibbs distribution on the manifold $\mathcal S_+^{n,p}$ equipped with two fundamental metrics. We have also provided several examples for which these sampling schemes can be numerically validated. 

\bibliographystyle{siamplain}  
\bibliography{reference}

\end{document}